\documentclass[11pt]{amsart}
\usepackage{amssymb,amsmath,epsfig,amsthm,color}
\date{}
\def\labelenumi{\theenumi)}
\def\C{{\mathbb C}}
\def\Z{{\mathbb Z}}

\def\R{{\mathbb R}}
\renewcommand{\geq}{\geqslant}
\renewcommand{\leq}{\leqslant}
\def\Tei{Teichm\"uller } 
\def\calT{\mathcal{T}}

\def\caht{c_{\scriptscriptstyle\mathrm{AHT}}}
\def\ll{\ell_\Delta}
\def\dw{\rho_{\Delta}}
\def\dl{\rho_{\mathcal L}}
\def\iso{\stackrel.=}
\def\comp{\asymp}
\def\co{\colon\thinspace}

\newtheorem{lemma}{Lemma}[section]
\newtheorem{theorem}[lemma]{Theorem}
\newtheorem{prop}[lemma]{Proposition}
\newtheorem{cor}[lemma]{Corollary}

\theoremstyle{definition}
\newtheorem{definition}[lemma]{Definition}

\theoremstyle{remark}
\newtheorem{remark}{Remark}[section]
\newtheorem{example}{Example}[section]

\begin{document}
\title{On the complexity of braids}
\author{Ivan Dynnikov}
\author{Bert Wiest}
\address{Dept.\ of Mech.\ and Math., Moscow State University,
Moscow 119992 GSP-2, Russia; LIFR MIIP, Bol'shoy Vlasyevsky 11,
Moscow 119002, Russia}
\email{dynnikov@mech.math.msu.su}
\address{IRMAR (UMR 6625 du CNRS), University of Rennes 1,
Campus de Beaulieu, 35042 Rennes cedex, France}
\email{\tt bertold.wiest@math.univ-rennes1.fr}
\begin{abstract}
We define a measure of ``complexity" of a braid which is
natural with respect to both an algebraic and a geometric
point of view. Algebraically, we modify the standard notion of the
length of a braid by
introducing generators $\Delta_{ij}$, which are Garside-like
half-twists involving strings $i$ through $j$, and
by counting powered generators $\Delta_{ij}^k$ as $\log(|k|+1)$
instead of simply $|k|$. The geometrical complexity is some
natural measure of the amount of distortion of the $n$ times
punctured disk caused by a homeomorphism. Our main result is that
the two notions of complexity are comparable.
This gives rise to a new combinatorial model for the \Tei
space of an $n+1$ times punctured sphere.
We also show how to recover a braid from its curve
diagram in polynomial time.
The key r\^ole in the proofs is played by
a technique introduced by Agol, Hass, and Thurston.
\end{abstract}
\maketitle

\tableofcontents

\section*{Introduction}

The $n$-string braid group $B_n$ can be defined in many different ways;
we shall mainly be interested in the following two:

Firstly, it has a finite presentation with generating set consisting
of Artin's half-twists $\sigma_1,\ldots,\sigma_{n-1}$. (The relations are
$\sigma_i\sigma_j=\sigma_j\sigma_i$ for pairs $(i,j)$ such that $|i-j|>1$,
and $\sigma_i\sigma_j\sigma_i=\sigma_j\sigma_i\sigma_j$ for $|i-j|=1$.)

Secondly, let us denote $D_n$ a closed disk in the complex plane, centered
on the origin, with $n$ punctures lined up on the real line. Then
we identify the braid group $B_n$ with the mapping class group of $D_n$
in the standard way \cite{Birman,DDRW}.
The idea of interpreting braids as isotopy classes of boundary-fixing
homeomorphisms of an $n$ times punctured disk $D_n$
is as old as the braid groups themselves.
Indeed, Artin~\cite{artin1,artin2}, who introduced the braid groups,
solved the word problem in $B_n$ by
showing that a braid $\beta\in B_n$ is uniquely characterized by
the images of generators of the fundamental group
of $D_n$ under the homeomorphism associated with $\beta$.

The aim of this paper is to clarify some aspects of the relation between
these two points of view on the braid group. Specifically, we show that
there is a notion of \emph{complexity} of a braid which is natural in both
frameworks, and which has an added advantage of allowing algorithmically
efficient transitions between the two perspectives.

The algebraic point of view is the following: as a generating set, we use
the set of all half-twists $\Delta_{ij}$ ($1\leqslant i < j \leqslant n$)
involving strands number $i$ though $j$. Now, the \emph{$\Delta$-length}
---our new algebraic notion of complexity---of a braid is simply the
shortest possible length of any word representing the braid, but with one
modification: a power $\Delta_{ij}^k$ ($k\in \Z$) shall not count as having
$\Delta$-length $|k|$, as would be usual, but as having $\Delta$-length
$\log_2(|k|+1)$.

The geometric measure of the ``complexity'' of a braid is as follows.
Let $E$ be a set of properly embedded disjoint simple arcs on the disk
$D_n$ separating all the punctures (we shall be using the set shown in
Fig.~\ref{F:exs}, but other diagrams would work just as well). Then a
\emph{curve diagram} of a braid $\beta$ is the image $\beta\cdot E$ of
$E$ under the homeomorphism $\beta$. A simple measure for the complexity
of a braid would be the number of intersections of $\beta\cdot E$ with
the real line (minimized over the isotopy class of $\beta\cdot E$).
However, since in a random braid word this quantity tends to grow
exponentially with the length, it is actually more natural not use this
quantity itself, but rather its logarithm as a measure of the complexity
of $\beta$.

The main result of this paper is that the two measures of complexity we
have just defined are comparable, in the sense that their ratio is
bounded from below and from above by positive constants depending only
on $n$.

In fact, our proof of this result also yields a new, algorithmically very
efficient way to calculate a \emph{canonical representative word} for any
element of the braid group $B_n$.

Until recently the most algorithmically efficient
treatments of the braid groups were based on purely
algebraic ideas that used a presentation of $B_n$
by generators and relations. The approach developed
by Garside~\cite{garside}, Thurston~\cite[Chapter~9]{epstein},
and Birman--Ko--Lee~\cite{bkl} yields an algorithm for
finding a canonical form of a braid $\beta\in B_n$ given initially
by a word in standard generators. If the input word is of length
$\ell$, the algorithm in \cite{bkl} requires at most $O(\ell^2n)$
operations. Note that Artin's original algorithm is actually very
inefficient in comparison, since the images of the generators
of $\pi_1(D_n)$ under the action of the given braid
may be of length $O(\exp(\mathrm{const}\cdot\ell))$.

The new algorithm described here, by contrast, is geometrical in flavor
and nevertheless efficient. The first crucial idea is that the curve
diagram of a braid can be computed efficiently. Indeed, it was pointed
out in~\cite{ybm}, \cite[Chapter~8]{DDRW} that a curve diagram
$D=\beta\cdot E$ can be naturally encoded by a $2n$--tuple of integers
$\eta(D)=(a_1,b_1,\ldots,a_n,b_n)\in\Z^{2n}$. Moreover, this vector can
be computed efficiently: the algorithm constructed there computes the
vector $\eta(\beta\cdot E)$ associated to a braid $\beta\in B_n$ of
length $\ell$ in time $O(\ell^2+n)$ (thus solving the braid recognition
problem in time $O(\ell^2)$, but without producing any kind of a
canonical word representing the given braid).

Having constructed the curve diagram of $\beta$, we are then looking for an
algorithm that, given the curve diagram $\beta\cdot E$ (or its associated
vector $\eta(\beta\cdot E)\in \Z^{2n}$) reconstructs a canonical word
representing the braid $\beta$. Moreover, we want this algorithm to be
efficient, and the output braid word to be about as short as possible, in
the sense that the $\Delta$-length of the output word surpasses the
minimal possible $\Delta$-length, among all representatives of $\beta$,
by only a linear factor. The centrepiece of this paper is an algorithm
which achieves just that.

The idea underlying our algorithm is simply to successively ``untangle''
the curve diagram $\beta\cdot E$. That is, given a curve diagram $D$,
one can act on it by a generator $\Delta_{ij}^k$ of the braid group so as
to simplify the diagram (or equivalently, such that the vector
$\eta(\Delta_{ij}^k\cdot D)$ is shorter than the vector $\eta(D)$, in
an appropriate metric). This process can then be repeated until the
diagram $E$ is reached. The braid word one has spelled out during the
untangling process is then a representative of $\beta^{-1}$. Such
untangling is always possible, but it is usually not unique, and the
difficulty is to do the untangling in an efficient manner. Our main
tool for doing so is a technique introduced by Agol, Hass, and Thurston
\cite{AHT}.

A note on history: the term ``curve diagram'' was introduced
in~\cite{FGRRW}, but the basic idea is much older: for instance, it is
very explicit in Mosher~\cite{MCGautom}, and indeed it is arguably
already present in Artin's original work. The fact that curve diagrams
are efficiently computable was popularized by one of us (I.D.) at the
{\it Braids\/} Colloquium in Toulouse, June 2000, and published later 
in~\cite{ybm}, \cite[Chapter~8]{DDRW}. It was also independently
discovered by Malyutin~\cite{mal} in slightly different settings.
We are not aware of any literature prior to that.
However, the fact that curve diagrams are determined
by their intersection numbers with a finite number of curves was
well-known before, see e.g.~\cite[Expos\'e 6]{FLP}.

The paper is organized as follows.
In Section \ref{S:motiv} we introduce two measures of complexity
of a braid, one geometric and one algebraic,  and formulate our main
result, that they are comparable.
In Section \ref{S:laminat} we introduce laminations---an important
instrument of our constructions. In Section \ref{S:AHT} we explain a
certain special case of Agol, Hass, and Thurston's algorithm for counting
the orbits of a collection of isometries of subintervals of $\Z$.
In Sections \ref{S:relax}, \ref{S:optimize} we show how this technique
can be combined with the idea of relaxing integral laminations in order
to prove the main theorem.
In Section \ref{S:twometrics} we rephrase the main result in terms of
quasi-isometries: we introduce two metrics on the braid group $B_n$
corresponding to the two measures of complexity, and prove that they are
quasi-isometric. In Section \ref{S:Teich} we prove that the metric space
constructed in the previous section embeds quasi-isometrically in the
Teichm\"uller space of the $n+1$ times punctured sphere, and is in fact
quasi-isometric to its so-called thick part.
In Section \ref{S:ordering} we show that our untangling
procedure provides an efficient algorithm for finding
$\sigma$-consistent braid word representatives.
In Section~\ref{S:algorithm} we discuss the complexity
of our algorithms. At the end of the paper, we discuss some perspectives 
for further research.


\section{Motivating example and statement of the result}\label{S:motiv}

As we told in the Introduction, we shall regard braids from $B_n$
as self-homeomorphisms of the punctured disk $D_n$, which are
viewed up to homotopy. The boundary $\partial D_n$ is supposed
to be fixed under all homeomorphisms that we consider.

We denote by $E$ the union of $n-1$ arcs in $D_n$ which
are shown in Fig.~\ref{F:exs} on the left. If $\beta$ is a braid,
then we let $\beta\cdot E$ be the union of arcs obtained from $E$
by the action of $\beta$, and we call this the \emph{curve diagram}
of $\beta$---this is only defined up to isotopies fixing the
boundary and the punctures. We recall, however, that by using such
an isotopy the curve diagram of any braid can be made \emph{tight}
with respect to the horizontal diameter of $D_n$ meaning that
each connected component of $\beta\cdot E$ and the real axis $\R$
are transverse to each other, and there are no puncture-free bigons
enclosed by them.
Each braid has a unique curve diagram which is tight with respect to
the horizontal diameter in the sense that any two such
diagrams are related by an isotopy of $D_n$ which
preserves the real axis. Details can be found e.g.\ in
\cite{FGRRW,DDRW}.

Throughout the paper, all curve diagrams we mention will be assumed
tight with respect to the axis unless otherwise specified.
We define the
\emph{norm}
of a curve diagram $D$ to be the
number of intersections of $D$ with the real axis:
$$\|D\|=\#(D\cap\R).$$

It is intuitively plausible that in order to create a very complicated curve
diagram, one needs a very long braid word. Equivalently, in order to
obtain the diagram $E$ by untangling a complicated curve diagram, one
needs to act on it by a long braid word. However, there is no simple
proportionality relation between length and complexity, as the following
example demonstrates.

\begin{example}\label{E:motivation}
Consider the following two braids: $\alpha=\sigma_2^{-1}\sigma_1$
and $\beta=\sigma_2\sigma_1$.
\begin{figure}[htb]
\centerline{\input{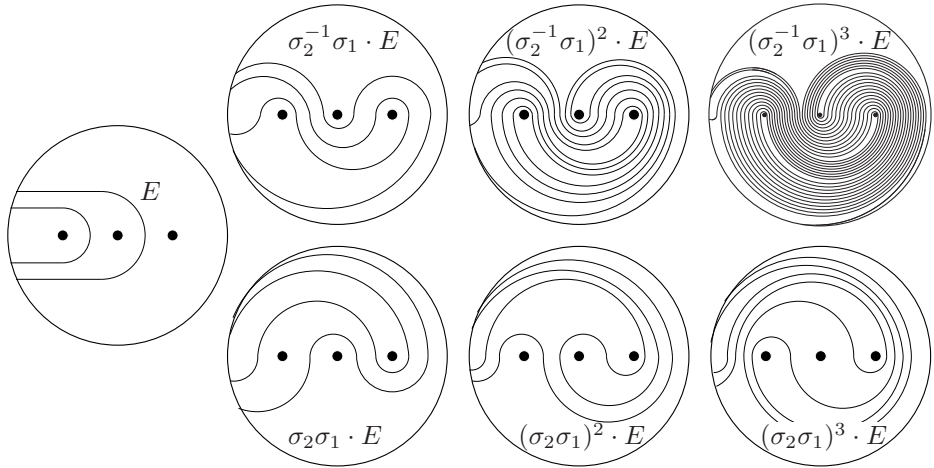}}
\caption{The diagram $E$ in $D_n$ consists of $n-1$ arcs, each intersecting
the real axis once
(shown here the case $n=3$).
The top row shows curve diagrams
of the braids $\alpha^k$, the bottom row of the braids
$\beta^k$.}\label{F:exs}
\end{figure}

The crucial observation now is that the sequence
$\|\alpha^k\cdot E\|$
grows
exponentially with $k$, whereas the sequence
$\|\beta^k\cdot E\|$
grows only linearly with $k$. Indeed, it is an exercise
to prove that $\|\alpha^k\cdot E\|=2(F_{k+2}-1)$,
where $F_0=1$, $F_1=1$, $F_2=2,\dots$
is the Fibonacci sequence. By contrast,
we have $\|\beta^k\cdot E\|=2\left[\frac{4k-1}3\right]+4$,
where $[x]$ stands for the integral part of $x$.
Thus,
$$\|\alpha^k\cdot E\|\sim\mathrm{const}^k,\qquad
\|\beta^k\cdot E\|\sim\mathrm{const}\cdot k.$$
Notice that both braids, $\alpha^k$ and $\beta^k$,
have the same length $2k$ in Artin's generators $\sigma_i$,
meaning that the shortest word representing any of them
has length $2k$.

The reason why there is a principal difference in the growth
of the complexity of curve diagrams $\alpha^k\cdot E$
and $\beta^k\cdot E$ is that $\alpha$ is a so-called
pseudo-Anosov braid, whereas $\beta^3=\Delta^2$, where
$\Delta$ is the Garside fundamental braid, i.e.\ the half-twist
of all strands at once.

Very roughly speaking, applying repeatedly the same twist $\Delta$
entangles the curve diagram much more slowly than applying first one
twist, then another one, then again a different one etc.
\end{example}

This example motivates the following
modification of the notion of braid length.
First, for
$1\leqslant i<j\leqslant n$ let $\Delta_{ij}$
denote the half-twist of strands $i$ through $j$:
\begin{equation}
\Delta_{ij}=(\sigma_i\ldots\sigma_{j-1})(\sigma_i\ldots\sigma_{j-2})
\ldots\sigma_i.
\label{E:delta}\end{equation}
This is a generating set of $B_n$,
which contains Artin's standard generators: $\sigma_i=\Delta_{i,i+1}$,
and the Garside fundamental braid: $\Delta=\Delta_{1n}$.

\begin{definition}
By the $\Delta$-\emph{length} of a word $w$ of the form
\begin{equation}\label{Delta-form}
w=\Delta_{i_1j_1}^{k_1}\ldots\Delta_{i_sj_s}^{k_s},
\end{equation}
where $k_t\ne0$ and
$\Delta_{i_t,j_t}\neq \Delta_{i_{t+1},j_{t+1}}$
for all $t$, we shall mean
$$\ll(w)=\sum_{i=1}^s\log_2(|k_i|+1).$$
For a braid $\beta\in B_n$ we define
$$
\ll(\beta)=\min\{\ll(w)\ | \ \hbox{the word }w\hbox{ represents }\beta \}.
$$
\end{definition}

Obviously, for any braid $\beta$, we have
$$\ll(\beta)\leqslant\ell(\beta),$$
where $\ell$ denotes the ordinary braid
length. Note that the $\Delta$-length is in general
not an integer.

\begin{definition}
We define the \emph{complexity} of a braid $\beta$ as
$$c(\beta)=\log_2\|\beta\cdot E\|-\log_2\|E\|.$$
\end{definition}

One of the main results of this paper is the following.

\begin{theorem}\label{T:main1}
The complexity and the $\Delta$-length of a braid
are comparable. More precisely,
there exist constants $C_1,C_2$ such that
the inequalities
\begin{equation}\label{Eq:main}
c(\beta)\leqslant C_1\cdot\ll(\beta)
\qquad\mbox{and}\qquad
\ll(\beta)\leqslant C_2\cdot n\cdot c(\beta)
\end{equation}
hold for any $\beta\in B_n$.
\end{theorem}

The constants $C_1,C_2$ will be given in (\ref{Eq:c1}), (\ref{Eq:c2})
below. We stress that they are independent of the number
of strands $n$.

Sections~\ref{S:laminat} through \ref{S:optimize} are devoted to the
proof of Theorem~\ref{T:main1}.


\section{Proof of the main result}

\subsection{Integral laminations}\label{S:laminat}

Curve diagrams that we consider belong to a wide and very natural class
of objects called laminations, which are defined without
any reference to the braid groups.

\begin{definition}
An \emph{integral lamination} in $D_n$ is a non-empty union $L$ of
finitely many
disjoint simple closed curves and properly embedded
arcs in $D_n$ such that
\begin{enumerate}
\item there are no bigons enclosed by the arcs of $L$
and $\partial D_n$ with no puncture inside;
\item no closed component of $L$ bounds a disk with no or
just one puncture inside;
\item the endpoints of arcs of $L$ are distinct from
$\R\cap\partial D_n$ and there is an equal number of them
above and below $\R$.
\end{enumerate}
\end{definition}

In particular, Conditions 1) and 2) of this definition imply that
all curves that can appear as connected components of a curve diagram
have the property that they get tangled under the action of some braids.
There is only one exception to this rule: a closed curve ``parallel'' to
the boundary of $D_n$ may be present in an integral lamination, even
though it does not get tangled under any self-homeomorphism of $D_n$.

Sometimes we may call integral laminations just laminations for simplicity.
The general notion of a lamination will not be needed here.

Integral laminations are regarded
modulo isotopy in $D_n$ fixing the boundary.
As in the case of curve diagrams, by such an isotopy any
integral lamination can be made tight with respect to the
axis $\R$.

In what follows, all integral laminations are assumed to be tight
with respect to $\R$ unless otherwise specified.
In some cases, for technical reasons, we shall consider
laminations that are ``almost tight'' with respect to $\R$.
Namely, some laminations $L$ that we consider
have exactly one puncture-free bigon enclosed
by $L$ and $\R$. However, laminations
are always assumed to be transverse to $\R$.

We do not distinguish between two laminations
$L_1$ and $L_2$ if $L_2$ is obtained from $L_1$
by an isotopy of $D_n$ preserving $\R$. In this case we
write $L_1=L_2$. If $L_1$ and $L_2$ are related by an
isotopy of $D_n$ not necessarily preserving $\R$, we write
$L_1\iso L_2$. Thus, if both $L_1$ and $L_2$ are tight with
respect to $\R$, then $L_1\iso L_2$ implies $L_1=L_2$.

The set of isotopy classes of laminations in $D_n$ will be denoted
by $\mathcal L_n$. The group $B_n$ acts on $\mathcal L_n$
in the natural way.

\begin{definition}\label{D:ccrelaxed}
A connected component of an integral lamination is
said to be \emph{relaxed} if it
intersects the axis $\R$ at most twice.
A lamination $L$ is said to be \emph{relaxed} if
all connected components of $L$ are relaxed.

A connected component of a lamination is called \emph{even}
(respectively, \emph{odd}), if it intersects the axis an even
(respectively, odd) number of times.
\end{definition}

If a lamination $L$ has the form $\beta\cdot E$, then each connected
component of $L$ is odd. Notice that if a lamination $L$ is very
complicated, then it may be quite difficult to decide if all its
components are odd, until we have untangled it.

We shall also make use of the following technical definition.

\begin{definition}
For an integral lamination $L$ in $D_n$ (not necessarily tight
with respect to $\R$), by a \emph{closure} of $L$
we shall mean the union $\overline L$ of pairwise disjoint
simple closed curves in the complex plane such that $\overline L\cap
D_n$ coincides with $L$, and $\overline L\setminus D_n$
is either empty or consists of arcs intersecting the axis $\R$ exactly
once, on the left of $D_n$. We shall view $\overline L$
up to a homeomorphism of the plane preserving the horizontal axis.
Clearly, $\overline L$ is unique up to such homeomorphisms.
\end{definition}

Curve diagrams as defined in Section~\ref{S:motiv}
are particular cases of laminations.
Our basic curve diagram $E$ is an example of a relaxed lamination.

The \emph{norm} of a lamination (not necessarily tight
with respect to $\R$) is defined in the same way as
for curve diagrams: $\|L\|=\#(L\cap\R)$. Notice:
$L_1\iso L_2$ does not necessarily imply $\|L_1\|=\|L_2\|$,
if at least one of the laminations $L_1$ or $L_2$ is
not tight with respect to the real axis.
However, as we have mentioned all laminations $L$
are assumed to be tight if not otherwise specified.

\begin{lemma}\label{L:one-Delta}
For any lamination $L$ in $D_n$, any $1\leqslant i<j\leqslant n$,
and $k\in\Z$, the following holds:
$$\|\Delta_{ij}^k\cdot L\|\leqslant(2|k|+1)\cdot\|L\|.$$
\end{lemma}

The proof, which is easy, will be left as an exercise to the reader.

\begin{proof}[Proof of the easy part of Theorem~\ref{T:main1}]
Let $\beta\in B_n$ be presented by the word~(\ref{Delta-form}).
For $t=0,\dots,s$,
put $\beta_t=\Delta_{i_{t+1}j_{t+1}}^{k_{t+1}}
\ldots\Delta_{i_sj_s}^{k_s}$,
so that $\beta_s=1$ and $\beta_0=\beta$.
Then Lemma~\ref{L:one-Delta} implies
$$\begin{aligned}
c(\beta)&=\log_2\|\beta\cdot E\|-\log_2\|E\|\\
&=\sum_{t=1}^s(\log_2\|\beta_{t-1}\cdot E\|-\log_2\|\beta_t\cdot E\|)\\
&=\sum_{t=1}^s(\log_2\|\Delta_{i_tj_t}^{k_t}\cdot(\beta_t\cdot
E)\|-\log_2\|\beta_t\cdot E\|)\\
&\leqslant\sum_{t=1}^s\log_2(2|k_t|+1)\\
&\leqslant\log_23\cdot\sum_{t=1}^s\log_2(|k_t|+1)
=\log_23\cdot\ll(w),
\end{aligned}$$
which gives the first inequality in~(\ref{Eq:main}) with
\begin{equation}\label{Eq:c1}
C_1=\log_23.
\end{equation}
\end{proof}

The difficult part of Theorem~\ref{T:main1}, i.e.~the
second inequality in~(\ref{Eq:main}), will be a consequence
of the following claim.

\begin{theorem}\label{T:main2}
For any integral lamination $L$ in $D_n$, there exists a
braid $\beta\in B_n$ such that $\beta\cdot L$ is relaxed
and the following holds
$$\ll(\beta)\leqslant9n\cdot\log_2\|L\|.$$
\end{theorem}

\begin{remark}
Note that if a lamination does not have the form of a curve
diagram, then there may be more than one braid untangling it.
\end{remark}

Thus the constant $C_2$ in~(\ref{Eq:main}) can be set to
\begin{equation}\label{Eq:c2}
C_2=9.
\end{equation}

Before embarking on the proof of the theorem, we make a remark
on how \emph{not} to prove this result. A naive approach could
be to try to proceed by analogy with the proof of the easy
part of Theorem~\ref{T:main1}, namely, to show that for a non-relaxed
integral lamination $L$ there exists a braid $\beta$ of the form
$\Delta_{ij}^k$ such that the following holds:
$$ \log_2\|\beta\cdot L\| \leqslant \log_2\|L\|-c\cdot\ll(\beta) $$
with some positive constant $c$ independent of $L$, and then apply
induction. Unfortunately, this does not work, because no such constant
$c$ exists. In other words, the function
$$\psi(L)=\inf_{\scriptscriptstyle\beta=\Delta_{ij}^k;\
\scriptscriptstyle\|\beta\cdot L\|<\|L\|}\
{\ll(\beta)}/(\log_2\|L\|-\log_2\|\beta\cdot L\|)$$
is unbounded. For instance, for $L=\Delta_{12}^N\Delta_{34}^N\cdot E$,
we have that $\psi(L)$ is comparable to $\log N/(\log2N-\log N)=
\log_2N$, and thus grows without bound as $N\rightarrow\infty$.
Moreover, changing the definition of the norm $\|{\cdot}\|$ in any naive
way does not help.

Instead of going this way, we shall apply a very
powerful technique due to Agol, Hass, and Thurston \cite{AHT}: we use a
certain complexity function which depends not only on the lamination,
but also on the history of the untangling process so far.


\subsection{The orbit-counting algorithm of Agol, Hass, and
Thurston}\label{S:AHT}

In this section, we give a brief account of Agol, Hass and Thurston's
technique, adapted to the special case which is relevant for our
purposes, and rephrased in a language which is more convenient for us.

By $[i,j]$ with $i,j\in\Z$ we denote the sequence
$$i,i+\epsilon,i+2\epsilon,\dots,j,$$
where $\epsilon=\pm1=\mathrm{sign}(j-i)$, and call such a sequence
an \emph{interval}. By the \emph{length} of the interval
$[i,j]$ we shall mean the number of elements in it, i.e., $|i-j|+1$.

For two intervals $[i,j]$ and $[k,l]$ of equal length, we define
the \emph{interval identification} $[i,j]\leftrightarrow[k,l]$
as the following symmetric relation in $\Z$:
$$i+p\epsilon_1\leftrightarrow k+p\epsilon_2,\qquad\mbox{for all }
p=0,\dots,|i-j|,$$
where $\epsilon_1=\mathrm{sign}(j-i)$ and $\epsilon_2=\mathrm{sign}(l-k)$.
(The interval identification $[i,j]\leftrightarrow[k,l]$ is not
distinguished from $[j,i]\leftrightarrow[l,k]$ and $[k,l]\leftrightarrow[i,j]$.)

\begin{definition}
An \emph{interval identification system} (\emph{IIS} for short)
$S$ is an interval $[1,N]$, $N>0$, endowed with a finite collection
of interval identifications
\begin{equation}\label{Eq:IIS}
[i_t,j_t]\leftrightarrow[k_t,l_t],\qquad t=1,\dots,r,
\end{equation}
within it: $i_t,j_t,k_t,l_t\in[1,N]$.
The number $N$ is called the \emph{norm} of $S$ and denoted $\|S\|$.
\end{definition}

For such an interval identification system $S$ we denote
by $\sim_S$ the equivalence relation in $[1,N]$
generated by the aggregate of
all identifications (\ref{Eq:IIS}).
The set $[1,N]/{\sim_S}$ of the equivalence classes will
be denoted by $\Omega_S$.

The Agol--Hass--Thurston algorithm that we are going to adapt
computes the size of $\Omega_S$ in time polynomial
in $(\log N,r)$. Note that all naive algorithms one may think
of immediately consume time linear in $N$, which is much worse.

The Agol--Hass--Thurston machinery was originally developed
for counting the number of connected components of a normal surface
given by its Haken coordinates. We shall apply it to simpler
geometrical objects and for quite a different purpose.

Namely, we shall consider the IISs that are naturally associated with
(the closure of) a lamination $L$ cut by a ray $(-\infty,a)\subset\R$.
The connected components of the cut lamination define an identification
relation between the intersection points $L\cap (-\infty,a)$.
The formal definition is as follows.

\begin{definition}\label{D:carries}
Let $\overline L$ be the closure of a lamination $L$ in $D_n$ which is
not necessarily tight with respect to the real axis.
Let $P_1,\dots,P_M$, $M=\|\overline L\|$ be the intersections
of $\overline L$ with the axis $\R$, numbered from left to right,
and let $S$ be an IIS with $N=\|S\|\leqslant M$.
We say that $L$ \emph{carries} $S$ if the following holds:
\begin{enumerate}
\item $i\leftrightarrow j$ holds in $S$ if and only if
the points $P_i$ and $P_j$ are connected by a segment of $\overline L$
not passing through a $P_t$ with $t\leqslant N$;
\item for any interval identification $[i_t,j_t]\leftrightarrow
[k_t,l_t]$ in $S$, there is a strip $R_t$ in the complex plane
bounded by the straight line segments $P_{i_t}P_{j_t},
P_{k_t}P_{l_t}\subset\R$ and two segments
$\gamma_t,\gamma_t'$ of $\overline L$
with $\partial\gamma_t=\{P_{i_t},P_{k_t}\}$,
$\partial \gamma_t'=\{P_{j_t},P_{l_t}\}$. (The arcs $\gamma_t$
and $\gamma_t'$ are allowed to coincide, in which case $R_t$
is just an arc.) The {\it bases} $P_{i_t}P_{j_t}$, $P_{k_t}P_{l_t}$
of the strip $R_t$ may overlap and even coincide; besides that,
$R_t$ must be embedded;
\item the strips $R_t$ are disjoint from each other
except at the axis $\R$.
\end{enumerate}
\end{definition}

Not every IIS is carried by a lamination.
For example, being carried by a lamination imposes the obvious
restriction that every integral point in the interval $[1,N]$
is involved in exactly two interval identifications, which
is not true in general.
Thus, what we consider is a particular case of the Agol--Hass--Thurston
algorithm.

For the rest of the paper, we shall never consider IISs or integral
laminations in isolation, but always an IIS $S$, carried by an
integral lamination $L$. Thus in our situation, it will be convenient
to use the geometrical language instead of the combinatorial one.
In particular, we shall refer to the elements of
$\Omega_S$ as connected components of $\overline L$ rather
than equivalence classes for $\sim_S$.

We shall assume that the whole picture of an integral
lamination is rescaled
so that the points $P_t$, $t=1,\dots,N$ coincide with
the integral points $1,\dots,N$ in the real axis.
Formally, by \emph{rescaling} we mean a homeomorphism
of the plane of the form $(x,y)\mapsto(\varphi(x),y)$,
where $\varphi$ is an increasing function. Clearly,
a rescaling preserves the combinatorial structure of a lamination.

We shall also speak of a \emph{strip} $[i,j]\leftrightarrow [k,l]$
instead of an interval identification $[i,j]\leftrightarrow
[k,l]$. Notice that such a strip carries a little more information
than the corresponding interval identification because in the
complex plane it can be attached ``from above'' and ``from below''
to the intervals $[i,j]$ and $[k,l]$. We shall always
assume that this above-below information is included
in the structure of the ISS $S$.

By the \emph{width} of a strip $R=([i,j]\leftrightarrow[k,l])$
we shall mean the number of connected components
of $L\cap(R\setminus \partial_0R)$, where $\partial_0R$
stands for the union of the bases of $R$.
Thus the width of $R$ is equal to the length of the bases
$[i,j]$, $[k,l]$ of $R$, i.e.\ $|i-j|+1$.
(Note that a strip of width one geometrically looks like
a strip of zero width.)

\begin{example}\label{E:IIS}
Shown in Fig.~\ref{F:transmex} is the closure of the lamination
$L=(\sigma_2^{-1}\sigma_1)^2\cdot E$ endowed with
the following interval identification
systems:
\begin{itemize}
\item[(a)]$N=26$,
$\{[1,6]\leftrightarrow[12,7]$, $[13,19]\leftrightarrow
[26,20]$, $[1,1]\leftrightarrow[26,26]$, $[2,2]\leftrightarrow[3,3]$,
$[4,14]\leftrightarrow[25,15]\}$;
\item[(b)]$N=25$,
$\{[1,6]\leftrightarrow[12,7]$, $[14,19]\leftrightarrow
[25,20]$, $[1,1]\leftrightarrow[13,13]$, $[2,2]\leftrightarrow[3,3]$,
$[4,14]\leftrightarrow[25,15]\}$;
\item[(c)]
$N=19$,
$\{[1,6]\leftrightarrow[12,7]$, $[4,9]\leftrightarrow[14,19]$,
$[1,1]\leftrightarrow[13,13]$, $[2,2]\leftrightarrow[3,3]$,
$[10,14]\leftrightarrow[19,15]\}$;
\item[(d)]
$N=14$,
$\{[1,6]\leftrightarrow[12,7]$, $[4,4]\leftrightarrow[14,14]$,
$[1,1]\leftrightarrow[13,13]$, $[2,2]\leftrightarrow[3,3]$,
$[5,9]\leftrightarrow[14,10]\}$.
\end{itemize}
In each picture, there are five strips in total,
two of them of width one, except in (d), where three strips are of width one.

\begin{figure}[htb]
\centerline{\input{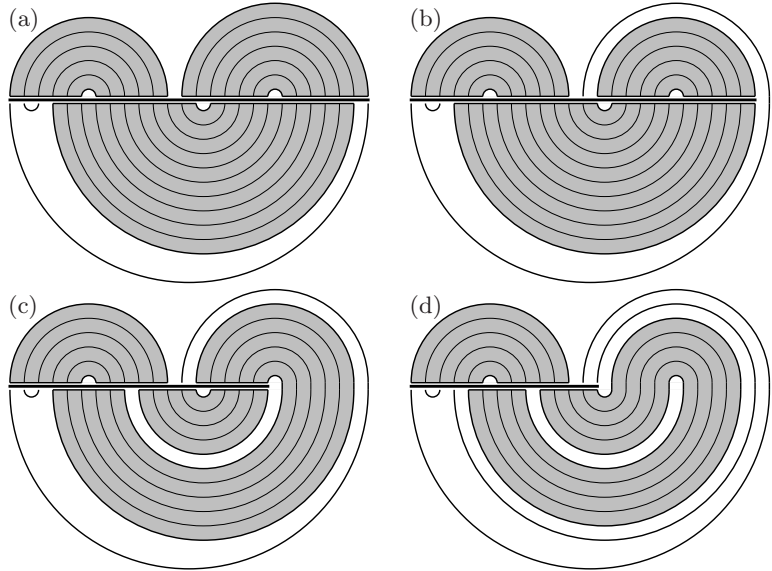}}
\caption{Four IISs carried by the lamination
$(\sigma_2^{-1}\sigma_1)^2\cdot E$.
The sequence (a)$\to$(b)$\to$(c)$\to$(d) is obtained from (a) by
successive transmissions.}
\label{F:transmex}
\end{figure}
In all these pictures, the strip $[1,6]\leftrightarrow[12,7]$ is attached
to both bases from above. The strip $[1,1]\leftrightarrow[13,13]$
in pictures (b)--(c) is attached to $\{1\}$ from below and to $\{13\}$
from above.
\end{example}

As this example demonstrates,
a lamination may carry many different interval
identification systems. If $L$ carries $S$ with given $\|S\|=N$, then
the relation $\sim_S$ is completely defined by $L$. This might suggest
that $S$ is uniquely defined by $L$ and $N$, but this
is not necessarily so. Indeed, the structure of $S$ assumes fixing
a collection of interval identifications, and therefore one can
genuinely change an IIS by replacing
an interval identification $[i,j]\leftrightarrow[k,l]$
of width at least two by two interval identifications
$[i,p]\leftrightarrow[k,q]$, $[p+\epsilon,j]\leftrightarrow[q+\epsilon',l]$,
where $p\in[i,j-\epsilon]$, $\epsilon=\mathrm{sign}(j-i)$,
$\epsilon'=\mathrm{sign}(l-k)$. Geometrically,
this means that some strip $R_t$ has been cut into two parallel strips.
This transformation truly changes the IIS while leaving both $N$ and
the underlying lamination $L$ invariant.

The idea of the Agol--Hass--Thurston's orbit counting algorithm is to
successively simplify an IIS by so-called \emph{transmissions}.
In a sense, this algorithm is a generalization of the well-known
Euclid algorithm for finding the greatest common divisor of integers.

In our specific situation it works as follows.
The input is an interval identification system $S$ which is carried by an
integral lamination $L$. We define a connected component counter, which
we set initially to zero.
At the ``rightmost point'' $N=\|S\|$ of the interval $[1,N]$ there are
exactly two strips
attached, one from below and the other from above.
Suppose $R_t=([i_t,N]\leftrightarrow[k_t,l_t])$
is the wider one of those two (or if they happen
to have the same width, let $R_t$ be either one of the two, no matter
which one).

\begin{figure}[b]
\centerline{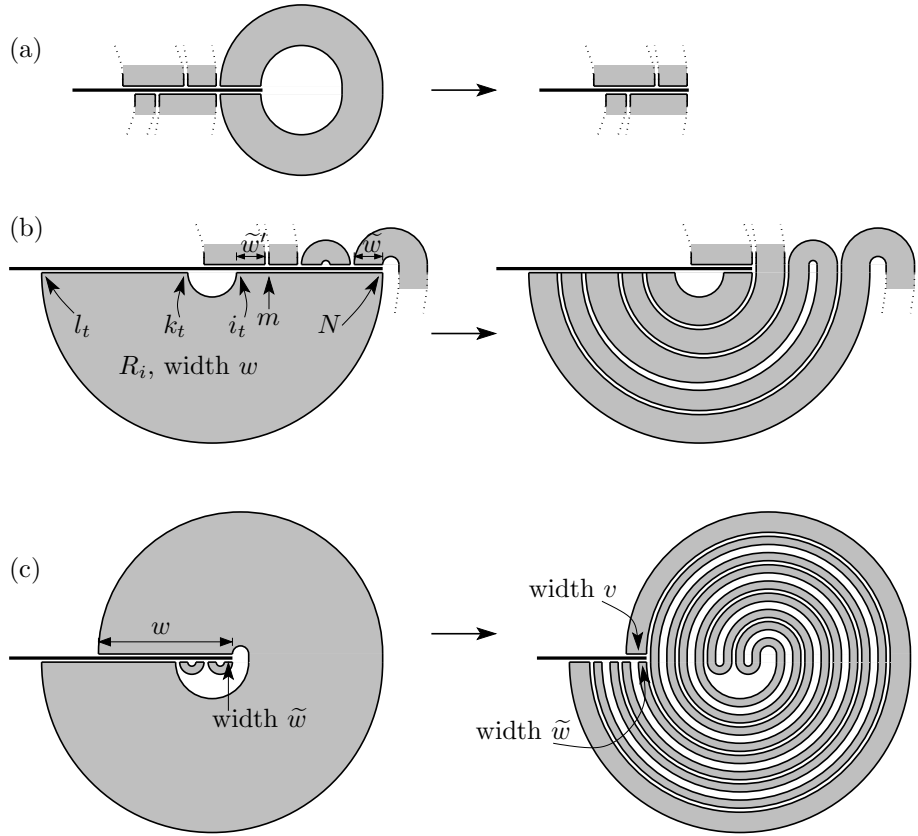}
\caption{(a) Removing an annulus; (b) a non-spiralling transmission;
(c) a twice-spiralling transmission}
\label{F:transm}
\end{figure}

It may even happen that those two strips coincide, which
means that $k_t=i_t$, $l_t=N$.
In this case, we advance the connected component counter
by $N-i_t+1$ and remove the strip $R_t$ from $S$. At the same time,
we replace $N$ by $i_t-1$. Intuitively, after this operation some
components of the lamination $L$ are no longer ``covered'' by any
strip of $S$. This operation is called
\emph{removing an annulus}. It is illustrated in Fig.~\ref{F:transm}(a).

Let us assume now that $l_t\ne N$. In this case we perform a
\emph{transmission}---a certain transformation which we are going to
explain next. Consider the set $X$ of all subintervals of $[i_t,N]$
that are bases of strips different from $R_t$. These strips will be
referred so as the \emph{denominators} of the transmission, whereas
the strip $R_t$ will be called the \emph{numerator}.
Let $m\in[i_t,N]$ be the left endpoint of the leftmost interval from $X$.
We cut the strip $R_t$ into a collection
of parallel strips so that the base $[i_t,N]$
is cut precisely into subintervals from $X$, and one more interval
$[i_t,m-1]$ provided that $m>i_t$. Then all pairs of strips whose bases
have just become matched are stuck together into longer strips,
and $N$ is set to $m-1$. This operation, which is called
transmission, is illustrated in Fig.~\ref{F:transm}(b), and
examples are given in Fig.~\ref{F:transmex}.

Thus, as a result of a transmission, the numerator and
all the denominators are replaced with strips that all,
except at most one, are obtained from the denominators by
attaching connected components of the cut numerator $R_t$.
If $m>i_t$, then there is one more strip of width $m-i_t+1$,
which we call the \emph{remainder} of the transmission.

There is one exception to the above rule: if the two bases
of $R_t$ overlap, i.e.\ we have $m-1=l_t>i_t$, then
we can perform the previously explained transmission
$d=\left[\frac{N-i_t+1}{N-m+1}\right]$ times at once. This is called a
$d$-\emph{times spiralling transmission}, and it
is illustrated in Fig.~\ref{F:transm}(c). The denominators and the
remainder of such a transmission are defined similarly to
those in the non-spiralling case.
In particular, the width of the remainder equals
$(N-i_t+1)-d(N-m+1)$.

It is obvious that under a transmission, the number of
elements in $\Omega_S$ does not change, the norm $N=\|S\|$ decreases, and the
number $r$ of interval identifications in $S$ does not increase.
Under an annulus removal operation, both $N$ and $r$ decrease,
and $|\Omega_S|$ decreases by the value added to the counter.
So, after finitely many operations described above, we end
up with an empty IIS, and then the connected component counter
indicates the number of connected components of the initial IIS.

It is not at all
obvious, however, that this procedure is efficient---in particular, that
only a relatively small number of transmissions is performed in the process.

\begin{definition}\label{D:AHT-compl}
The \emph{AHT-complexity} $\caht(S)$ of a nontrivial IIS
$$S=\{[i_t,j_t]\leftrightarrow[k_t,l_t]\;|\;t=1,\dots,r\}$$
is
$$\caht(S)=r+\sum_{t=1}^r\log_2(|i_t-j_t|+1)-\frac12\log_2\widetilde w,$$
where $\widetilde w$ is the width of the narrower strip attached
to $N=\|S\|$.
If $S$ is the trivial ISS, we put $\caht(S)=0$.
\end{definition}

\begin{remark}
In this definition,
the summand $-\frac12\log_2\widetilde w$ looks quite artificial, and
it was not present in the original definition by Agol, Hass, and
Thurston. Introducing it allows to prove a better estimation for the
simplification effect of a transmission in our specific case of
laminations. Note that what we subtract is just
one half of one of the summands in the preceding sum, so
we count the logarithm of the width of one selected strip with weight
one half, whereas all the logarithms of
other widths are counted with weight one.
\end{remark}

\begin{lemma}\label{L:AHTred}
(a) Performing a non-spiralling transmission on $S$
reduces $\caht(S)$ by at least one.

(b) Performing a $d$-times spiralling transmission on $S$ reduces
$\caht(S)$ by at least $\log_2(d+1)$.
\end{lemma}

\begin{proof}
(a) The value $\widetilde w$ in Definition~\ref{D:AHT-compl}
is the width of the ``rightmost'' denominator of the transmission
to be applied to $S$. Let $\widetilde w'$ be the width of the rightmost
denominator in the next transmission step.
In addition, let $w$ be the width of
the numerator of the transmission.

Suppose that the remainder of the transmission is not trivial.
Then its width is exactly $\widetilde w'$, and we have
$w\geqslant\widetilde w+\widetilde w'$. The transmission
causes the following change of the AHT-complexity:
$$\caht^{\mathrm{old}}-\caht^{\mathrm{new}}=
\log_2w-\frac12\log_2\widetilde w-\frac12\log_2
\widetilde w'=\frac12\log_2\frac{w^2}{\widetilde w\widetilde w'}
\geqslant1,$$
which follows from the fact that $(a+b)^2\geqslant4ab$ for all $a,b>0$.

If the remainder is trivial, then the new IIS has a smaller
number of interval identifications, so we have
$$\caht^{\mathrm{old}}-\caht^{\mathrm{new}}\geqslant1+
\log_2w-\frac12\log_2\widetilde w+\frac12\log_2\widetilde w'
\geqslant1.$$

(b) Let $w,\widetilde w,\widetilde w'$ be as before. Suppose that the
remainder is nontrivial and let $v$ be its width. Then the
two strips attached to the rightmost point of the IIS obtained after
the transmission have widths $\widetilde w$ and $v$.

If $\widetilde w
\leqslant v$, then $\widetilde w'=\widetilde w$, and we have
$$\caht^{\mathrm{old}}-\caht^{\mathrm{new}}=
\log_2w-\log_2v\geqslant\log_2(d+1),$$
because $w/v\geqslant d+1$.

If $\widetilde w\geqslant v$, then $\widetilde w'=v$, and we have
$$\caht^{\mathrm{old}}-\caht^{\mathrm{new}}=
\log_2w-\frac12\log_2\widetilde w-\frac12\log_2v
\geqslant\frac12\log_2\frac{(d\widetilde w+v)^2}{\widetilde wv}\geqslant
\log_2(d+1),$$
since $w\geqslant d\widetilde w+v$ and
$(na+b)^2/ab>(n+1)^2$ for all $a\geqslant b$, $n\geqslant1$.

Finally, if there is no remainder, then we have
$$\caht^{\mathrm{old}}-\caht^{\mathrm{new}}=1+
\log_2w-\frac12\log_2\widetilde w+\frac12\log_2\widetilde w'\geqslant
\log_2d+1\geqslant\log_2(d+1),$$
since $w\geqslant d\widetilde w$.
\end{proof}


\subsection{Relaxing integral laminations}\label{S:relax}

In this section we prove the following claim, which is a ``weaker version''
of Theorem~\ref{T:main2}.

\begin{theorem}\label{T:main2a}
For any integral lamination $L$ in $D_n$, there exists a
braid $\beta\in B_n$ such that $\beta\cdot L$ is either relaxed
or contains a relaxed even component, and
the following holds
$$\ll(\beta)\leqslant8n^2\cdot(\log_2\|L\|+1).$$
\end{theorem}

\begin{remark}
For a lamination $L$ that is not a curve diagram,
Theorem~\ref{T:main2a} asserts that by a braid of $\Delta$-length
$O(\log_2\|L\|)$ we can ``partially'' untangle $L$ so that an even
component is revealed. It is actually possible to untangle such an $L$
completely by a braid of the indicated
$\Delta$-length. The proof of this fact
requires more technical details,
which we prefer to postpone until the next section.
\end{remark}

The proofs of Theorems~\ref{T:main2} and \ref{T:main2a} follow essentially
the same scheme, but the argument for Theorem~\ref{T:main2a} is more
``straightforward". So, in this section, we explain the main principle
that allows to prove an inequality of the form
$\ll(\beta)\leqslant\mathrm{const}(n)\cdot\log_2\|L\|$,
whereas the next section contains details that allow to
make the $\mathrm{const}(n)$ grow as slowly as $O(n)$.

The basic idea is this: we think of our IIS as being
made of a rigid horizontal line and a number of rubber-rectangles
attached to it, and after each transmission we allow the picture to
``relax''.

More rigorously, by \emph{relaxing} a lamination $L$ we mean applying
a braid $\beta$ so that the lamination gets simpler,
i.e., so as to have $\|\beta\cdot L\|<\|L\|$.
For any curve diagram distinct from $E$,
there may be many braids of the form $\Delta_{ij}^k$ that relax it,
and it is very easy to find at least one of them. However,
recursively applying relaxing braids of the form $\Delta_{ij}^k$
in a naive way until a relaxed lamination is reached
may result in an untangling braid word of length $O(\|L\|)$,
because it can be only guaranteed
that each relaxation reduces the norm of $L$ by at least some
additive constant. The use of the AHT algorithm allows
to make a choice of a relaxation at each step of the untangling process
so that the $\Delta$-length of the untangling braid word is
of order $O(\log\|L\|)$.

Our algorithm then works as follows: for a given integral lamination $L$,
we construct an IIS $S_0$ such that $L$ carries $S_0$ and
$\|S_0\|=\|\overline{L}\|$.
Then the construction of the previous section yields a sequence
\begin{equation}
S_0\stackrel{d_1}\longmapsto S_1\stackrel{d_2}\longmapsto
S_2\stackrel{d_3}\longmapsto \ldots\stackrel{d_p}\longmapsto
S_p=\varnothing,
\end{equation}
where by $\stackrel d\longmapsto$ with $d\geqslant2$ we denote a
$d$-times spiralling transmission, by $\stackrel1\longmapsto$ a
once-spiralling transmission or a non-spiralling transmission,
and by $\stackrel0\longmapsto$ the removal of an annulus.
Lemma~\ref{L:AHTred} implies
\begin{equation}
\sum_{i=1}^p\log_2(d_i+1)\leqslant\caht(S).
\end{equation}
Moreover, for any $i=1,\dots,p$ the IIS $S_i$ is still carried by the
lamination $L$.

We put $L_0=L$, and subsequently find laminations $L_1,L_2,\dots,L_p$
such that, for any $i=1,\dots,p$, the following holds:
\begin{enumerate}
\item $L_i$ carries $S_i$;
\item either $L_i=L_{i-1}$ or
$L_i$ is obtained from $L_{i-1}$ by a relaxation,
$L_i\iso\beta_i\cdot L_{i-1}$;
\item $L_i$ is the simplest lamination satisfying 1) and 2)
(in the sense that it has the minimal norm).
\end{enumerate}
It is not required here that all the laminations $L_i$ are tight with
respect to $\R$. We remark, however, that they will be not far from
being tight, and most of them will actually be tight. Indeed, the reader
who just wishes to understand the principle of the algorithm may safely
forget about non-tight laminations.

If $L$ has only odd components, then we end up with
$$L_p=\beta\cdot L=E,$$
where $\beta=\beta_p\beta_{p-1}\ldots\beta_1$.
The desired relation between $\ll(\beta)$ and $\|L\|$ is then
obtained by estimating $\caht(L)$ and $\ll(\beta_i)$ for all $i$.

\begin{example}
Consider again the closure $\overline L$ of the curve
diagram $L=(\sigma_2^{-1}\sigma_1)^2
\cdot E$. Fig.~\ref{F:transmex}(a) shows the corresponding IIS $S_0$,
and from the picture, one can guess the general rule for choosing $S_0$.

The systems of strips shown in Fig.~\ref{F:transmex}(b), (c), (d)
correspond to the lamination-IIS pairs $(L,S_1)$, $(L,S_2)$,
$(L,S_3)$, respectively. One can easily see that no braid will
simplify $L$ if we require that the obtained lamination still carry
$S_1$ or $S_2$. This is because the strips in Fig.~\ref{F:transmex}(b), (c)
are relaxed, i.e.~embedded in the plane ``in the optimal way" with respect
to the number of intersections with the real axis.

\begin{figure}[htb]
\centerline{\input{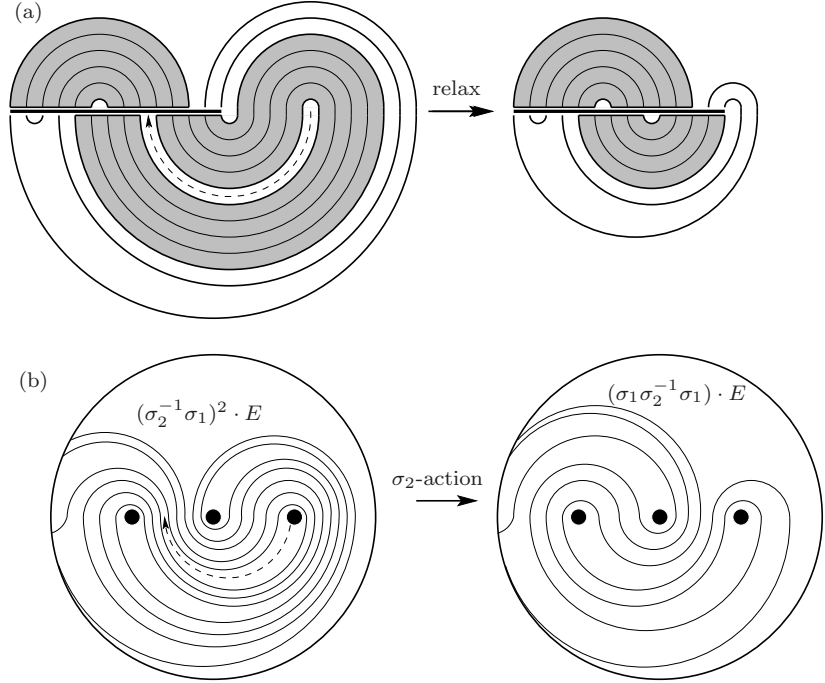}}
\caption{Part (a) of this figure is a continuation of Fig.~\ref{F:transmex}
(the left hand side is the same as Fig.~\ref{F:transmex}(d)).
It gives an example of a ``relaxation''. Part (b) shows the corresponding
relaxation of the curve diagram of $(\sigma_1\sigma_2^{-1})^2$.}
\label{F:relaxex}
\end{figure}

In Fig.~\ref{F:transmex}(d), by contrast, one of the strips has two
``unnecessary" intersections with $\R$, which can be cancelled if we apply
the half-twist $\sigma_2$. This is illustrated in Fig.~\ref{F:relaxex}.
Thus, in this example, we have $L_1=L_2=L$, $L_3=\sigma_2\cdot L$.
\end{example}

Now we give a formal description of the construction. Recall
that we assume the whole picture of the lamination
to be rescaled so that the closure $\overline L$
intersects the real axis in the points
$1,2,\dots,N=\|\overline L\|$.

First of all, we need to define $S_0$. In order to do so, we cut the
lamination $\overline L$ by the whole real axis, thus obtaining an
identification relation $\leftrightarrow$ on the interval $[1,N]$,
where $N=\|\overline L\|$. Then we collect each maximal
family of parallel arcs of the cut lamination
into a single strip of $S_0$. In other words, we choose $S_0$
carried by $L$ so that $\|S_0\|=\|\overline L\|$ and
$S_0$ has the minimal possible
number of interval identifications (strips).

\begin{lemma}\label{L:inicomp}
We have
$$\caht(S_0)\leqslant(2n-1)(\log_2\|L\|+1).$$
\end{lemma}

\begin{proof}
The width of any strip of any IIS carried by $L$ is obviously no larger
than the norm of $L$ (unless $L$ is already relaxed).
Therefore, it suffices to show that there are
at most $(2n-1)$ strips in $S_0$.
This can be done in numerous ways, e.g., as follows.

On $S^2=\C\cup\infty$, take a foliation $\mathcal F$ with singularities,
such that:
\begin{enumerate}
\item all connected components of $\overline L$ are leaves of $\mathcal F$;
\item $\mathcal F$ is transverse to the segment $I=(1,N)\subset\R$
except at the punctures;
\item all the singularities of $\mathcal F$ are simple
(see Fig.~\ref{F:foliation}) and the number of them
is minimal provided that there is a singularity of
type~1 at $\infty$.
\end{enumerate}
\begin{figure}[ht]
$$\begin{array}{ccc}
\epsfig{file=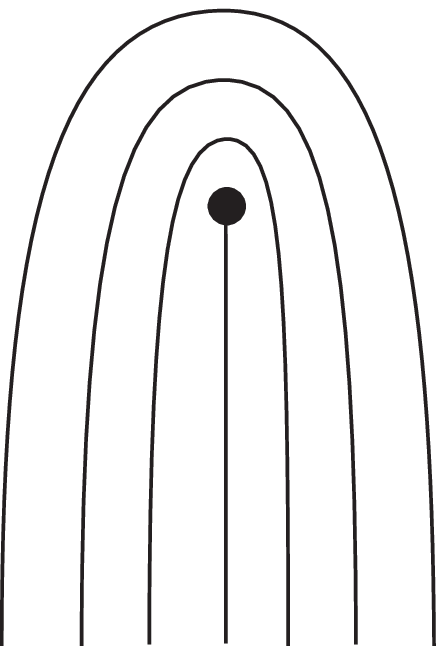,height=70pt}&\hskip2cm&
\epsfig{file=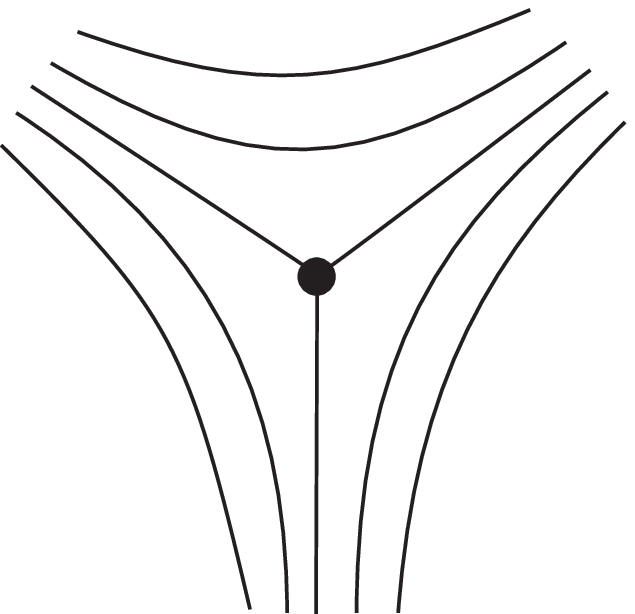,height=70pt}\\
\text{Type 1}&&\text{Type 2}
\end{array}$$
\caption{Singularities of $\mathcal F$}\label{F:foliation}
\end{figure}

We denote by $k_1$ and $k_2$ the number of singularities of type~1 and 2
shown in Fig.~\ref{F:foliation}.
Then the singularities of type~1
may occur only at punctures, at the leftmost point of $D_n$,
and at infinity, so we have $k_1\leqslant n+2$.
Moreover, we have $k_2=k_1-4$, which follows from the Poincar\'e--Hopf
formula $\chi(S^2)=\frac{1}{2}k_1-\frac{1}{2}k_2$.

If a region complementary to $I\cup\overline L$, and not containing
$\infty$, contains $s\geqslant1$ singularities of type~2 then it is
adjacent to $s+2\leqslant3s$ strips of $(L,S)$ lying on the same
side of $\R$. If the outermost region, in which $\infty$ lies,
contains $s\geqslant 0$ singularities of type~2, then it is
adjacent to $s+3\leqslant 3s+3$ strips.
The sides of strips that are not part of the
boundaries of the just mentioned regions are
in one-to-one correspondence with singularities
of type~1 lying in $I$. Thus, for the number $r$ of strips, we have
$$2r\leqslant 3(k_1-4)+3+(k_1-1)\leqslant 4n-2,$$
which completes the proof of the lemma.
\end{proof}

Now we proceed with describing the untangling process.
Each pair $(L,S)$, where $L$ is a lamination carrying the IIS $S$,
defines a collection of strips. This motivates the following notation:

\begin{definition}\label{D:ssrelaxed}
(a) A \emph{strip system} is a pair $(L,S)$, where $L$ is a lamination
and $S$ is an IIS carried by $L$.
(In a strip system $(L,S)$ the lamination $L$ is not necessarily
assumed to be tight with respect to $\R$.)

(b) For a strip $R_t$ of a strip system, we call the number of
connected components in $R_t\setminus\R$ the \emph{length} of $R_t$.

(c) A strip system $(L,S)$ is said to be \emph{relaxed} if
all strips in it are of length $\leqslant2$ and all the connected
components of $L$ that are not covered by strips
are relaxed.
\end{definition}

Clearly, a strip system is not relaxed if and only if there is
a deformation of the complex plane that preserves the bases of
the strips and makes at least one of them shorter. The idea of
such a deformation is to make the strips tight with respect to
the axis $\R$. However,
it may be impossible to find such a deformation
if we require punctures to be fixed
during the deformation. So, in a sense,
some punctures provide an obstruction to the relaxation.
All the obstructing punctures are located on the right
of the right base of the remainder.

\begin{definition}
Let $\alpha$ be a semicircular arc in the complex plane such
that $\partial\alpha=\{P,Q\}\subset\R$, where $P$
is a puncture and $Q$ is not.
By \emph{sliding the puncture} $P$ along $\alpha$
we mean a homeomorphism $\varphi$ of the complex plane such that
\begin{enumerate}
\item $\varphi$ is identical outside a small neighbourhood
$U$ of $\alpha$;
\item $U$ does not contain any other puncture except $P$;
\item $\varphi$ takes $P$ to $Q$.
\end{enumerate}
\end{definition}

Viewed up to rescaling, each sliding a puncture
operation gives rise to a braid.
The crucial observation now is that the corresponding braid
can be decomposed into two or fewer $\Delta$s:
$$\Delta_{ij}^\epsilon\Delta_{i,j\pm1}^{-\epsilon},$$
where $\epsilon=\pm1$ and, by definition, $\Delta_{ii}=1$.
Such a braid is called semicircular in \cite{Wi}.

\begin{figure}[ht]
\centerline{\epsfig{file=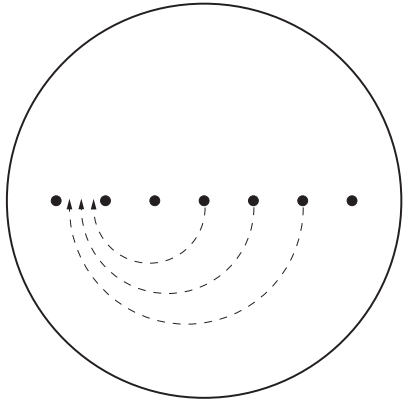,height=150pt}}
\caption{The sliding of these three punctures is
represented by $\Delta_{26}\Delta_{23}^{-1}$
and by $\Delta_{56}^{-1}\Delta_{26}$}\label{F:sliding}
\end{figure}

It is also important to note that in some cases, we can slide
a few punctures simultaneously by applying a braid of $\Delta$-length
$\leqslant2$. This occurs if we slide the punctures
along ``parallel'' arcs, and there are no other punctures
between the moving ones and between their destinations,
see Fig.~\ref{F:sliding}. The corresponding braid can be
represented in the form $\Delta_{ij}^\epsilon\Delta_{ij'}^{-\epsilon}$
and in the form $\Delta_{ij}^\epsilon\Delta_{i'j}^{-\epsilon}$.

We are now ready to describe the untangling process completely.
For a given a lamination $L$, we start with
finding the initial IIS $S_0$.
We set $L_0=L$.

\begin{figure}[ht]
\centerline{\epsfig{file=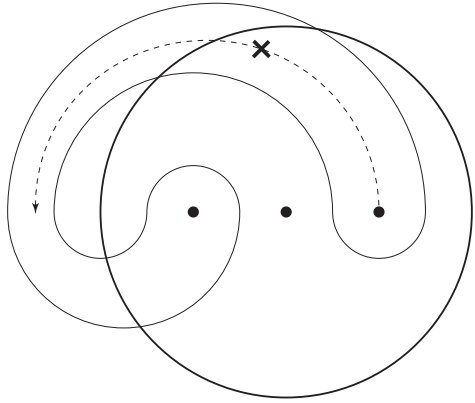,height=150pt}}
\caption{Punctures should not be slid out of $D_n$}\label{F:forbid1}
\end{figure}

Then, for each $i=1,2,\dots$ we do the following.
First, we examine the strip system $(L_{i-1},S_{i-1})$,
which is relaxed by construction. One of the following
situations must occur:
\begin{enumerate}
\item
$S_{i-1}=\varnothing$, which means that $L_{i-1}$ is
relaxed. In this case, we terminate the procedure.
\item
An annulus removal operation applies to $S_{i-1}$.
The strip to be removed contains only relaxed components of $L_{i-1}$.
If at least one of them is even, we terminate the procedure.
If all of them are odd, we remove the annulus, set $L_i=L_{i-1}$,
$\beta_i=1$, and proceed as before.
One can actually see that at this point the lamination is
already untangled, so, after removing a few annuli the process will
be terminated.
\item
A transmission $S_{i-1}\stackrel{d_i}\longmapsto S_i$ applies
to $S_{i-1}$. If the strip system $(L_{i-1},S_i)$ is still relaxed,
we set $L_i=L_{i-1}$, $\beta_i=1$ and proceed as before. If not,
we slide the punctures that obstruct the relaxation, along
arcs parallel to
the arcs of $L$ until they reach the bases of some strips.
After that we deform $L_{i-1}$, keeping the new positions of
punctures fixed, so as to reduce the number of intersections
with $\R$ on the right of $N_i=\|S_i\|$ as much as possible.
(In most cases this just means to make the lamination tight
with respect to $\R$ for the new positions of punctures. However,
an example of a situation where this is \emph{not} the case
is given in Fig.~\ref{F:ab}.) This replaces $L_{i-1}$
by $L_i\iso\beta_i\cdot L_{i-1}$, where the braid $\beta_i$ is
obtained by combining all the slidings. The strip system $(L_i,S_i)$
is now relaxed, and we proceed as before. However,
there is an exception: if $L_{i-1}$ contains an even non-closed
component, then we may not be able to slide some obstructing
puncture so as to let the strip system get
relaxed, see Fig.~\ref{F:forbid1}. In this case,
we apply the slidings until we have that even component of $L$
relaxed in the sense of Definition~\ref{D:ccrelaxed}
(although the strip system is not yet relaxed in
the sense of Definition~\ref{D:ssrelaxed})
and terminate the procedure.
One last note: sometimes there may be more than one way to slide
an obstructing puncture. If so, we choose the way that
allows to move the puncture farther to the right.
\end{enumerate}

The key ingredient of the proof of Theorem \ref{T:main2a} is the
following bound on the lengths of the braids $\beta_i$.

\begin{lemma}\label{L:relaxlen}
For all $i=1,2,\dots,p$ we have
\begin{equation}\label{E:relaxlen}
\ll(\beta_i)\leqslant 4n\cdot\log_2(d_i+1).
\end{equation}
\end{lemma}

\begin{proof}
The assertion is nontrivial only in the transmission case,
$d_i\geqslant1$. We consider the non-spiralling case ($d_i=1$) first.
In Fig.~\ref{F:onerelaxn} all possible types of obstructing punctures
are indicated.
In each case, we need to slide such a puncture at most twice.
Since the number of punctures to be slid is not larger than $n$,
we have in this case:
$$\ll(\beta_i)\leqslant 2\cdot 2\cdot n=4n\cdot\log_2(d_i+1).$$

\begin{figure}[htb]
\centerline{\input{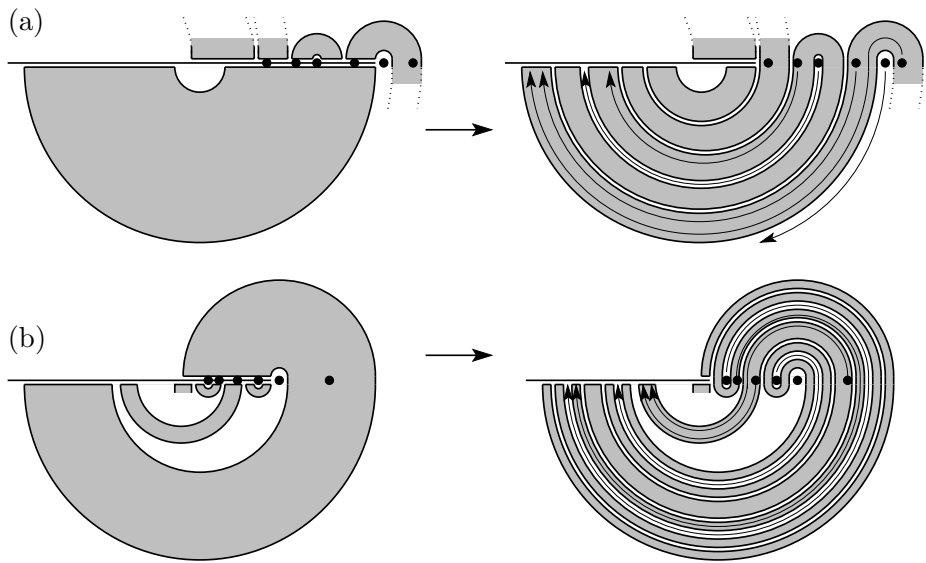}}
\caption{How to relax after a transmission, in the non-spiralling case.}
\label{F:onerelaxn}
\end{figure}

\begin{figure}[htb]
\centerline{\input{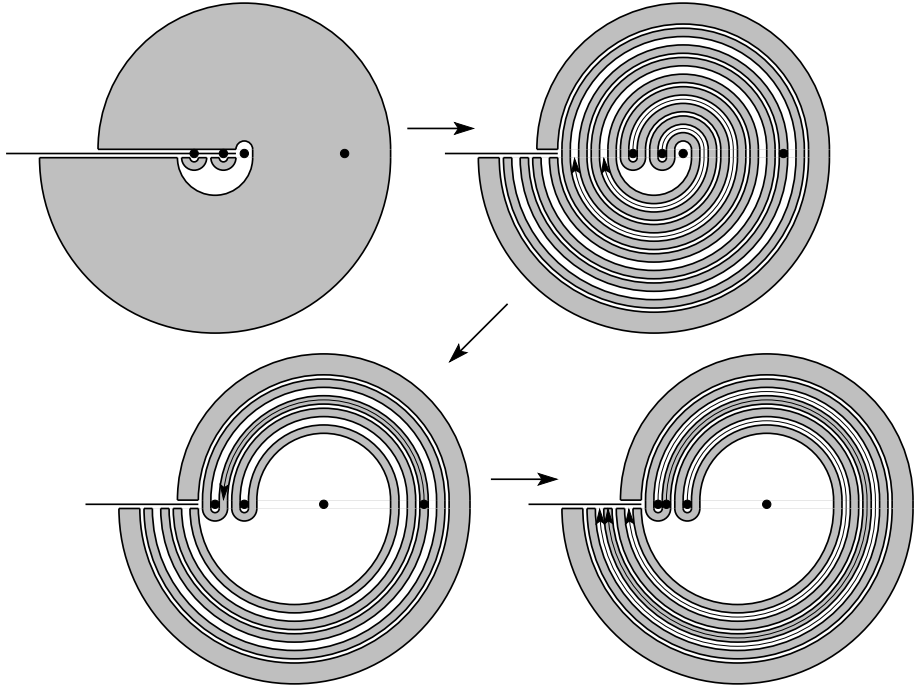}}
\caption{How to relax after a transmission, in the $d$-times spiralling case.}
\label{F:onerelaxn2}
\end{figure}

Now we turn to the case of a spiralling transmission
$S_{i-1}\stackrel{d_i}\longmapsto S_i$. Fig.~\ref{F:onerelaxn2}
shows how the obstructing punctures should
be slid. There are necessarily some punctures that must
be slid $2d$ times. If there are no more obstructing punctures
and no punctures at the left base of the numerator,
the whole spiral can be untwisted by a braid of the form
$\Delta_{ij}^{\pm 2d}$. There is, however,
a complication, if some obstructing punctures are at
a smaller distance than $2d$ from their final destination:
during the untwisting process, more and more
punctures need to be ``picked up''.
For instance, in Fig.~\ref{F:onerelaxn2}
we have a twice-spiralling relaxation, but after the first two half
turns of three punctures, a fourth puncture gets picked up.

Let there be obstructing punctures at distances
$k_1,k_2,\dots,k_q$ from their final destination,
$2d=k_1>k_2>\ldots>k_q\geqslant0$,
where we regard all punctures within the left base of
the numerator also as obstructing (even though they may
be at distance zero from their destination).
Without loss of generality we may assume that there is
exactly one puncture in the center of the spiral.

In order to deliver all obstructing punctures to their
destinations, we first apply the $(k_1-k_2-1)$st power of
the half-twist involving the ``farthermost'' punctures and
the central puncture.
Now we have to pick up the punctures that are at distance
$k_2$ from their destination. We do this by sliding them
back one step. Now they are in a row with the punctures from the first group
and all of them are at distance $k_2+1$ from
their destination. We apply the $(k_2-k_3)$th power
of the half-twist involving all these punctures and the central
puncture, then pick up the next portion of punctures, and so forth.

The total number of punctures picked up during the untwisting
is at most $n-2$,
and we have $q\leqslant n-1$. Picking up each puncture ``costs''
us two $\Delta$s.
Thus in this case the $\Delta$-length of $\beta_i$ is bounded
from above by
$$2(n-2)+\sum_{i=1}^q\log_2(k_i-k_{i+1})\leqslant
2(n-2)+(n-1)\log_2(2d)<4n\cdot\log_2(d+1),$$
where we set $k_{q+1}=-1$.
We leave it to the reader to treat the case when some even component
of $L$ gets relaxed during the untwisting of the spiral.
\end{proof}

\begin{proof}[Proof of theorem~\ref{T:main2a}]
This is now an easy corollary to Lemmas~\ref{L:inicomp}, \ref{L:AHTred},
and \ref{L:relaxlen}:
we start with an IIS $S_0$, carried by the lamination $L$, of
AHT-complexity at most $(2n-1)(\log_2\|L\|+1)$. Then we apply
alternatingly transmission moves and relaxation moves. During the $i$th
transmission, the AHT-complexity gets reduced by \emph{at least}
$\log_2(d_i+1)$, whereas in the subsequent relaxation move a braid of
$\Delta$-length \emph{at most} $4n\cdot\log_2(d_i+1)$ is applied.
Since the AHT-complexity cannot take negative values, the cycle must
stop before a braid word of $\Delta$-length $4n\cdot (2n-1) \cdot
(\log_2\|L\|+1)$ has been spelt out.
\end{proof}

\begin{remark}
The untangling process described in this section has been implemented as
a {\tt maple}-worksheet by Michel Bonnefont and Erwan Hillion. Their
program, which draws pictures of both the curve diagrams and the
interval identification systems, is freely available \cite{BonHil}.
\end{remark}


\subsection{Optimizing the untangling procedure}\label{S:optimize}
In this section we modify the arguments of the previous section
so as to obtain a proof of Theorem~\ref{T:main2}. We use a
very similar construction, but modify the definitions of
$S_0$, $\beta_i$, and $L_i$. In order to distinguish from the
previous constructions, we add a prime in the notation:
$S_0'$, $\beta_i'$, $L_i'$. Instead of Lemmas~\ref{L:inicomp} and
\ref{L:relaxlen}, we shall get the following estimates:
\begin{align}
\caht(S_0')&<3n\cdot\log_2\|L\|-3n,\label{E:newinicomp}\\
\ll(\beta_i')&\leqslant3\cdot\log_2(d_i'+1),\quad
i=1,\dots,p'\label{E:newrelaxlen}\\
\ll(\beta_{p'+1}')&\leqslant2n.\label{E:lastbeta}
\end{align}
which, together with Lemma~\ref{L:AHTred} imply Theorem~\ref{T:main2}.
Thus the achievement is to get rid of the factor $n$ in
the estimation (\ref{E:newrelaxlen}), which is the
counterpart of (\ref{E:relaxlen}), at the expense
of enlarging the constant in~(\ref{E:newinicomp})
and getting a more involved construction.

The reason for the factor $n$ appearing in (\ref{E:relaxlen}) is that
we don't know how many punctures we need to slide at each relaxation
step, and we estimate the number very roughly by $n$. The idea now
is to move almost all those punctures at once, using the trick
indicated in Fig.~\ref{F:sliding}. To this end,
we must make sure that there are no punctures in between the destination
points, so that the moved punctures do not get shuffled with the others.

Once the new untangling process is described, it is easy, though
tiresome, to verify that it works and relations (\ref{E:newinicomp}),
(\ref{E:newrelaxlen}), (\ref{E:lastbeta}) hold. We skip some details of
this checking, which contains not much new compared with the previous
section. What we do in detail is describing the new
rules for relaxing.

First we recall that laminations
and their closures that we consider
are forbidden to pass through the leftmost point of the disk $D_n$.
The reader might have noticed that the role of this point in
our figures is similar to those of punctures. Now it will
become even more similar. We call this point the \emph{false puncture}
and mark it by $*$ in the figures.

During the untangling process, we shall treat the
false puncture almost in the same as a ``true'' one.
Namely, we consider the closure $\overline L$ of the initial
lamination $L$ as an ordinary lamination in an $(n+1)$-punctured
disk $D_{n+1}$ whose punctures are the same as before plus
the false puncture.

Denote by $\iota$ the inclusion $B_n\rightarrow B_{n+1}$ given by
$\iota(\sigma_i)=\sigma_{i+1}$. At the $i$th step of the untangling process,
the relaxing braid $\beta_i'$ will be, in general, a braid
from $B_{n+1}$. However, the resulting braid $\beta_{p'+1}'\ldots\beta_1'$
will lie in $\iota(B_n):\beta_{p'+1}'\ldots\beta_1'=\iota(\beta)$.
This is achieved by organizing the untangling process so that
\def\theenumi{(\roman{enumi})}
\def\labelenumi{\theenumi}
\begin{enumerate}
\item\label{condition1}
a true puncture is never slid below the false one;
\item
the false puncture is never slid; the transmission-relaxation
procedure is terminated as soon as both bases of the numerator
of the transmission to be applied are on the left of $*$,
or we get $S_{p'}'=\varnothing$;
\item
once a true puncture has been moved to the left
of $*$, it stays untouched until the final step,
when all true punctures that have been
slid to the left of $*$ are slid towards the right of $*$ along
arcs in the upper half-plane;
the additional braid $\beta_{p'+1}'$ does this job.
\end{enumerate}

It is not hard to show (using Lemma~\ref{L:sigma-free} below) that
\begin{equation}\label{E:sum}
\ll(\beta)\leqslant\sum_{i=1}^{p'+1}\ll(\beta_i').
\end{equation}

A base of a strip will be called an \emph{A-base} if the strip
approaches it from above, and a \emph{B-base} otherwise
(`A' stands for `above' and `B' for `below').
To each strip, we associate its \emph{type} that can be either
AA, AB, BA, or BB depending on the types of the bases:
the first letter indicates the type of the left base, and the second
of the right one. If the bases of the strip coincide, it can
be thought of as an AB- or BA-strip, this does not matter.

A strip system $(L,S)$ is said to be \emph{almost relaxed}
if the length of all its BB-strips is not larger than three,
and for all the other strips not larger than two.
As before, $L$ is not assumed to be tight with respect
to the axis, but all puncture-free bigons enclosed
by $L$ and $\R$ must be on the right of $\|S\|$.

We define the new untangling procedure $\overline
L=L_0'\mapsto L_1'\mapsto\dots$
so as to comply with
the following rules:
\begin{enumerate}
\setcounter{enumi}3
\item
$S_0'\stackrel{d_1'}\longmapsto S_1'\stackrel{d_2'}\longmapsto
S_2'\stackrel{d_3'}\longmapsto \ldots\stackrel{d_{p'}'}\longmapsto S_{p'}'$
is a sequence of transmissions and annulus removal operations;
\item\label{condition5}
for any $i=0,\dots,p'$, the strip system $(L_i',S_i')$ is almost
relaxed;
\item
for any $i=1,\dots,p'$ we have $L_i'\iso\beta_i'\cdot L_{i-1}'$
with some $\beta_i'\in B_{n+1}$;
\item\label{condition7}
for any $i=0,\dots,p'$, the interiors of all strips of $(L_i',S_i')$
and their A-bases are free of punctures; the false puncture $*$
is not contained in any (A- or B-) base of a strip.
\end{enumerate}

The IIS $S_0$ from the previous section (which, we recall, has at
most $2n-1$ strips) does not in general satisfy
Condition~\ref{condition7}. This is because some punctures
may sit on the A-bases of strips, and a base
of a strip may contain $*$. We resolve this
by cutting those strips into a few parallel ones.
This results in enlarging the number of strips by at most $n+2$,
and one can show that the number of strips will be enlarged
exactly by $n+2$ only if it was strictly smaller than $2n-1$
before cutting. So, the number $r$ of strips in $S_0'$ is at most $3n$.

One now obtains~(\ref{E:newinicomp}) as follows.
Let $w_1,\dots,w_r$ be the widths of the strips of $S_0'$.
If $L$ is not relaxed, then $r\geqslant4$. We also have:
$$\qquad\|L\|\geqslant \frac r2, \qquad\sum_{i=1}^rw_i\leqslant\|L\|.$$
This implies
$$\begin{aligned}
\caht(S_0')&\leqslant r+\sum_{i=1}^r\log_2w_i
=r+\log_2\left(\prod_{i=1}^rw_i\right)\\
&\leqslant r+\log_2\left(\frac{\sum_{i=1}^rw_i}r\right)^r\\
&\leqslant r+r(\log_2\|L\|-\log_2r)\\
&=3n\cdot\log_2\|L\|-(3n-r)\log_2\|L\|-r\log_2\frac r2\\
&\leqslant3n\cdot\log_2\|L\|-3n\cdot\log_2\frac r2
\leqslant3n\cdot\log_2\|L\|-3n.
\end{aligned}$$

Provided that Conditions \ref{condition5}, \ref{condition7}
above are satisfied up to $i=k-1$, we shall explain
how to define $\beta_k'$.
For simplicity, we will assume that, during the untangling
process, no two punctures become immediate neighbours so
that the lamination does not traverse the interval between
them. One can easily show that this is not a loss of generality,
since such two punctures can be treated as a single one.

The notion of obstructing puncture was defined somewhat loosely
in the previous section. Now we make it more precise.
To this end, consider the strip system $(L_{k-1}',S_{k-1}')$
and the transmission $S_{k-1}'\stackrel{d_k'}\longmapsto S_k'$.

\begin{definition}
An arc $\alpha\subset L_{k-1}'$ will be said to be \emph{essential}
if it satisfies the following conditions:
$\alpha$ lies in the lower half-plane, and we have $\partial\alpha\subset\R$;
the left endpoint of $\alpha$ is located in the right base of the
numerator of the transmission $S_{k-1}'\stackrel{d_k'}\longmapsto S_k'$,
but not in the right base of the remainder.
\end{definition}

By definition, an essential arc is contained in
the numerator or in a denominator of the transmission
$S_{k-1}'\stackrel{d_k'}\longmapsto S_k'$.
We also remark that it has both its endpoints to the right of $*$.
Our untangling process will be organized as follows: at each step,
the essential arcs form a family of parallel, concentric semicircles;
in particular, it makes sense to talk about an outermost essential arc.
The strips of $(L_{k-1}',S_k')$ that require relaxation
after the transmission will be
exactly those that contain an essential arc of $(L_{k-1}',S_{k-1}')$.
The relaxation is achieved by ``pushing all essential arcs across the
real line". So, by \emph{obstructing punctures} we shall mean those
punctures that are located between the endpoints of the outermost
essential arc.

Depending on the type of the numerator of the transmission,
the following cases are possible:

\begin{figure}[htb]
\centerline{\epsfbox{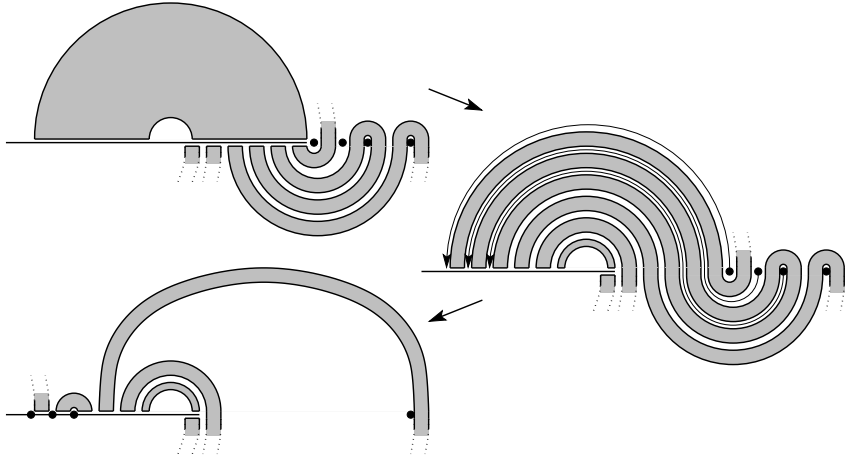}}
\caption{Case where the numerator is of type AA}\label{F:aa}
\end{figure}

\smallskip\noindent{\it Case AA.}
In this case BB-denominators of length one with both bases participating in
the transmission cannot occur. Indeed, between those bases there
must be a puncture, which contradicts to the requirement that
all A-bases are free of punctures. Thus, any length one
BB-denominator has one of its bases further to the left. Such a
denominator gives rise to a length two AB- or BA-strip,
which does not need to be simplified.

All the other denominators are of AB type and length two,
or BB type and length three.
The obstructing punctures should be slid along arcs parallel to
the denominators toward the right
base of the numerator (by one $\Delta$), and then
along the numerator toward the left base (two more $\Delta$s),
\begin{figure}[ht]
\centerline{\epsfbox{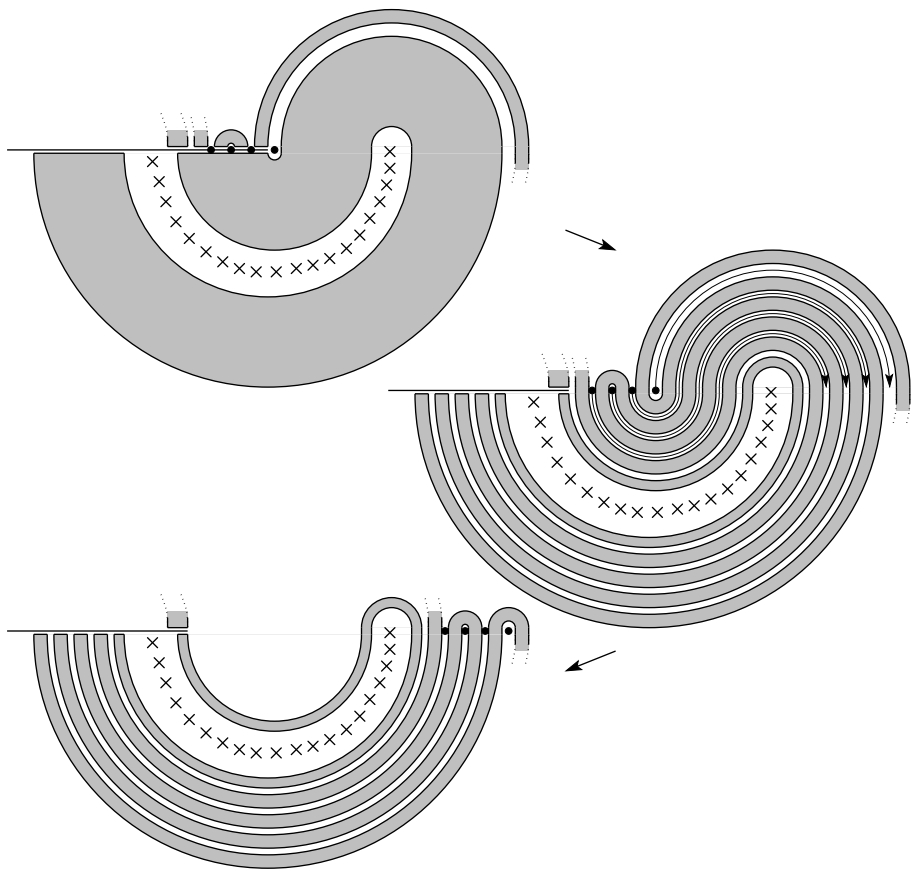}}
\caption{Case where the numerator is of type BB, length three}\label{F:bb}
\end{figure}
\begin{figure}[ht]
\centerline{\epsfbox{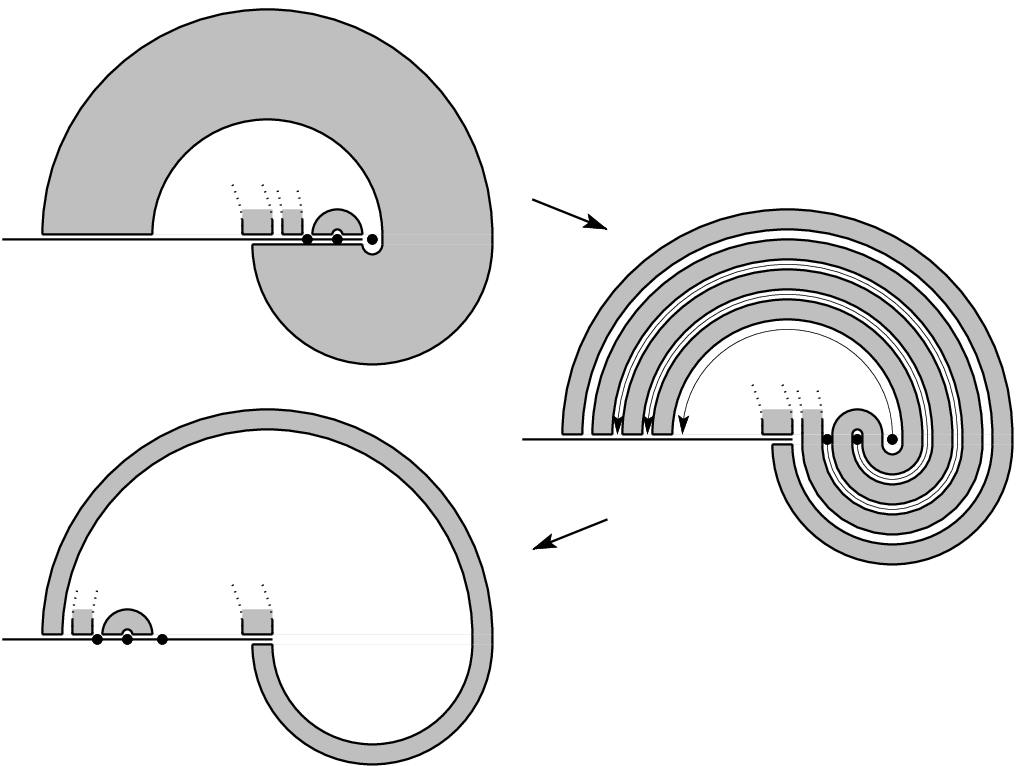}}
\caption{Case where the numerator is of type AB, non-spiralling}\label{F:ab}
\end{figure}
see Fig.~\ref{F:aa}. We make just one exception to this
rule: if all denominators are of
type BB, then an obstructing puncture positioned on the
immediate right of the right base of the numerator (if there is one)
does not participate in the second sliding. This is not important
for the moment but will be in the proof of Lemma~\ref{L:sigma-consistency}
below. In this case, we have $\ll(\beta_k')\leqslant3$.

\smallskip\noindent{\it Case BB, length one.}
No relaxation is needed
at this point, since every strip that is created during the transmission
is of AB type and length two, or of BB type and length three.

\smallskip\noindent{\it Case BB, length three.} The obstructing
punctures may be inside the right base of the numerator and
on the immediate right of that base. They are
slid twice along the numerator to the right, see Fig.~\ref{F:bb}.
We have $\ll(\beta_k')\leqslant3$.

\smallskip\noindent{\it Case AB, non-spiralling.}
The obstructing punctures, which are inside
and on the immediate right of
the B-base of the numerator, are slid twice along the numerator,
see Fig.~\ref{F:ab}. Again, we have $\ll(\beta_k')\leqslant3$.

\begin{figure}[hbt]
\centerline{\epsfbox{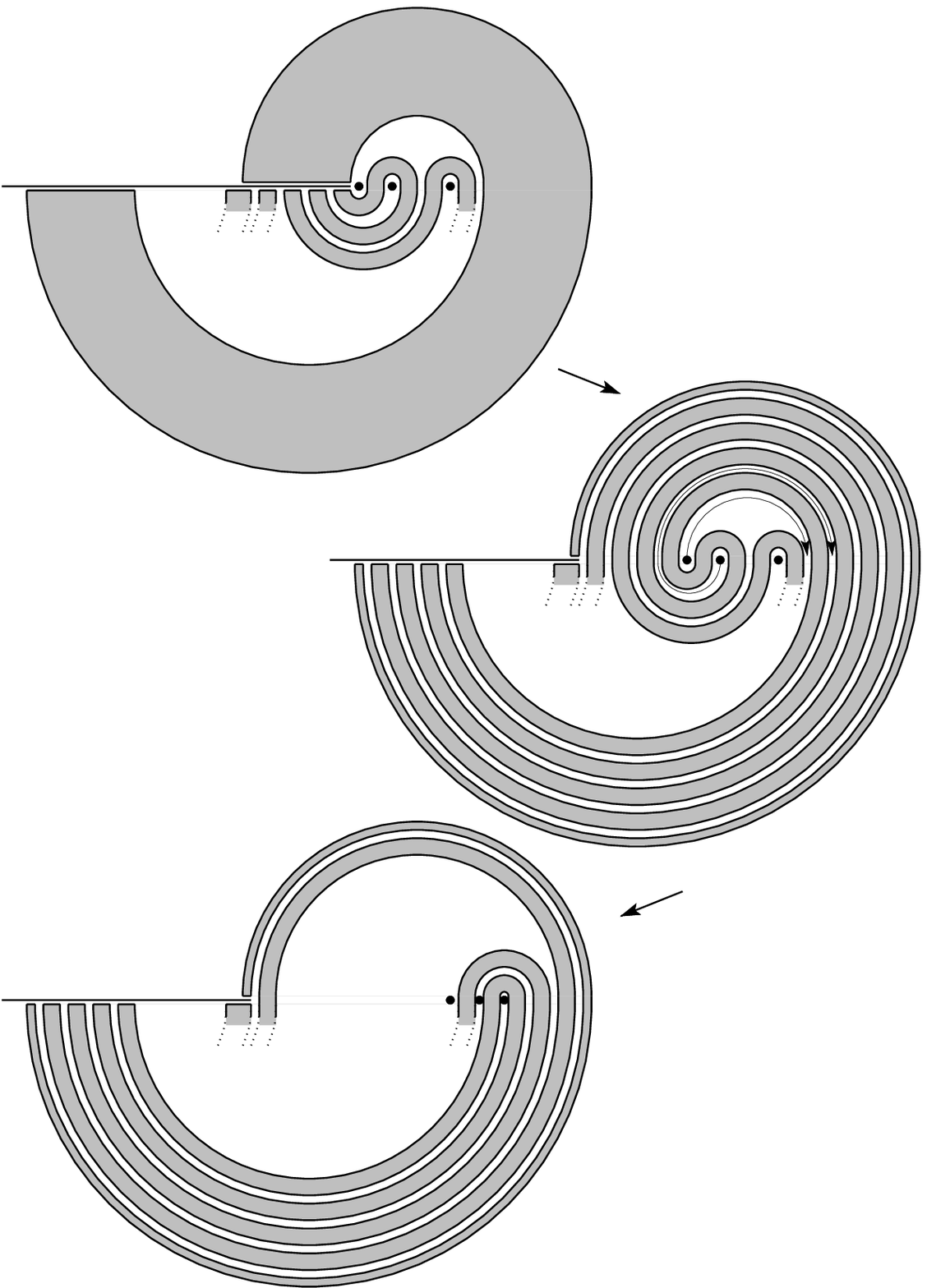}}
\caption{Case where the numerator is of type BA, non-spiralling}\label{F:ba}
\end{figure}

\smallskip\noindent{\it Case BA, non-spiralling.} All the denominators
are of BB-type. Those of length one must have the other base further
to the left. After the transmission,
they give rise to BB-strips of length three, which don't need to be
relaxed for the moment. The denominators of length three give
rise to strips of length five or seven, see Fig.~\ref{F:ba}.
The obstructing punctures are first slid along arcs parallel
to essential ones, and then once along the numerator.
As always, we have $\ll(\beta_k')\leqslant3$.

\smallskip\noindent{\it Spiralling case.}
The difficulty with punctures that need to be picked up (see the
previous section) now disappears, because there are no punctures
in the interior of any strips. Thus we can simply apply
$\Delta_{ij}^{2d_k'}$ in the BA-case and $\Delta_{ij}^{-2d_k'}$
in the AB-case, where the half-twist
$\Delta_{ij}$ involves the punctures inside the spiral.
We have $\ll(\beta_k')=\log_2(2d_k')<2\log_2(d_k'+1)$.
Note that in the BA-spiralling case all the denominators are of type
BB and of length three. After the relaxation, the strips they
give rise to are also of length three.

\smallskip
The transmission-relaxation process is terminated once
the transmission ``cutting edge'' has arrived at $*$.
Thus, during the process, whenever an obstructing puncture
is slid along an arc in the lower half-plane, the arc will be
above the outermost essential arc, and, therefore, on the right
of $*$. This guarantees that Condition \ref{condition1} holds.

As a result of the process, some true punctures
have moved to the left of $*$. At the very end of the relaxation
process, we slide them back along
arcs in the upper half-plane so
as to get the simplest possible lamination. This
yields a braid $\beta_{p'+1}'$ of $\Delta$-length
at most $2n$, because at most $n$ punctures need to be slid
(actually, it has $\Delta$-length at most $n$, but
even $8n$ would be good enough for our purposes).

\begin{figure}[ht]
\centerline{\raisebox{66.5pt}
{\hbox to 0pt{\hskip 28.5pt\Large\bf*\hss}}
\epsfig{file=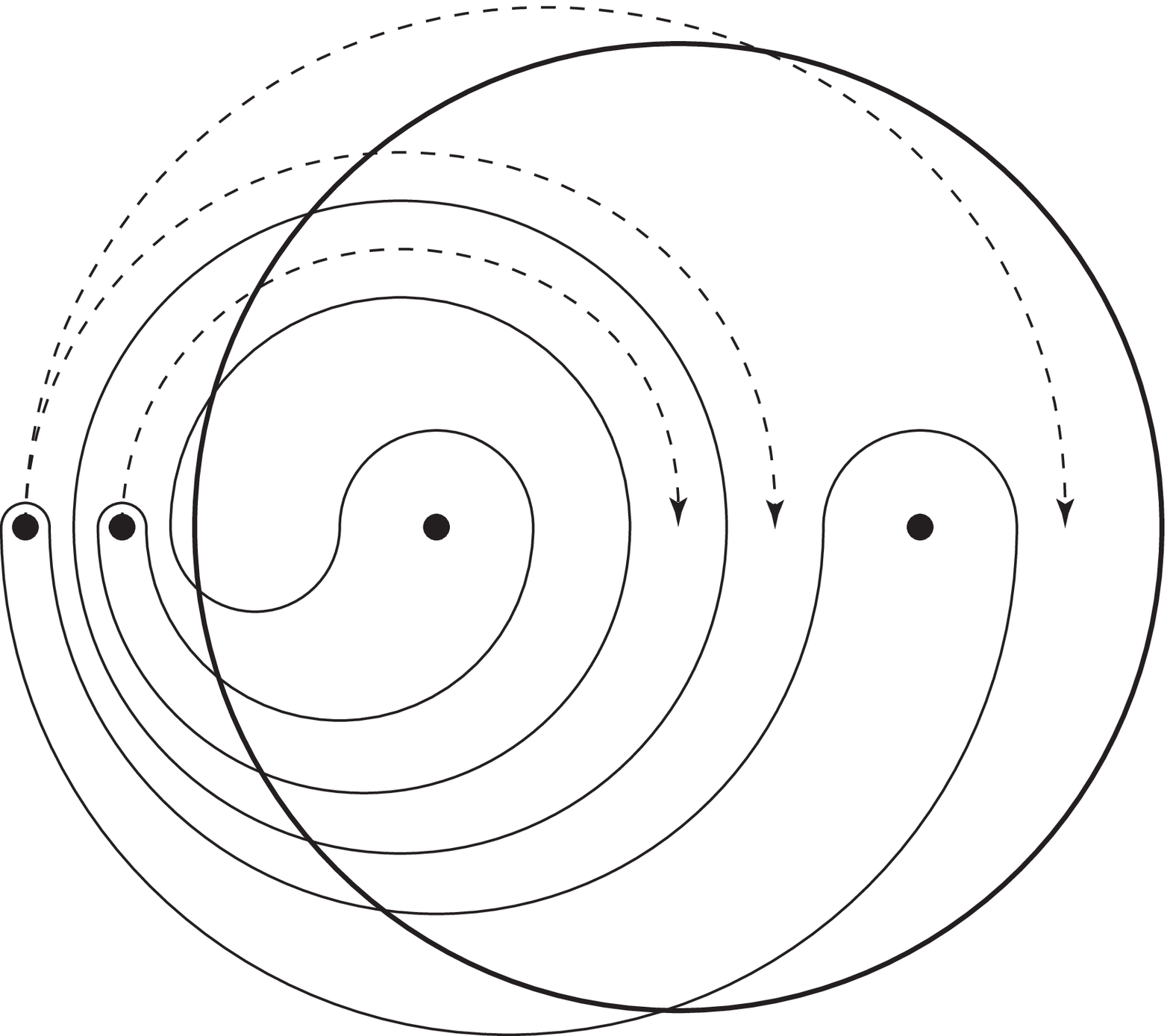,height=150pt}}
\caption{The braid $\beta_i'$ returns true punctures
into the disk. In the situation shown here, there are two possibilities for
the leftmost puncture, and just one possibility for the next
one}\label{F:lastbeta}
\end{figure}
We remark that if the original lamination $L$ contains even
components, then there may be no canonical choice for $\beta_{p'+1}'$,
see Fig.~\ref{F:lastbeta}.

This completes our description of the optimized untangling procedure,
and thus the proof of Theorem~\ref{T:main2}.

\begin{remark}
In each braid $\beta_i\in B_{n+1}$ just described, there
is one strand that corresponds to the false puncture $*$.
By removing this strand, we obtain a braid $\beta_i''\in B_n$,
and we have
$$\beta=\beta_{p'+1}''\ldots\beta_2''\beta_1''.$$
Consequently applying braids $\beta_1'',\beta_2'',\dots$ to $L$ may,
in principle, eventually result in a more complicated lamination than 
the original one. This may occur when $\beta_i''$ corresponds to 
sliding a puncture to the left of $*$. So, though the braid $\beta$ 
untangles the initial lamination $L$, the sequence
$$\|L\|,\|\beta_1''\cdot L\|,\|\beta_2''\beta_1''\cdot L\|,\ldots$$
may not be monotonically decreasing. \end{remark}


\section{Applications}
\subsection{Two equivalent metrics on $B_n$}\label{S:twometrics}

The $\Delta$-length of a braid defined in Section~\ref{S:motiv}
possesses the following obvious properties:
$$\ll(\beta)=\ll(\beta^{-1})\qquad\mbox{and}\qquad
\ll(\beta_1\beta_2)\leqslant\ll(\beta_1)+\ll(\beta_2).$$
This means that the function
\begin{equation}
\dw\colon\thinspace B_n\times B_n \to \R, \ \dw(\beta_1,\beta_2)=
\ll(\beta_1^{-1}\beta_2),
\end{equation}
is a distance on $B_n$.

The analogue statement for the complexity function is not true:
in general, $c(\beta)\neq c(\beta^{-1})$. However, it is
true for the following modified complexity function. Put
\begin{equation}
\widetilde c(\beta)=\sup_{L\in\mathcal L_n}
\bigl|\log_2\|\beta\cdot L\|-\log_2\|L\|\bigr|.
\end{equation}
By definition, for all $\beta\in B_n$, we have
$c(\beta)\leqslant\widetilde c(\beta)$ and
$$\begin{aligned}
\widetilde c(\beta)&=\sup_{L\in\mathcal L_n}
\bigl|\log_2\|\beta\cdot(\beta^{-1}\cdot L)\|-\log_2\|\beta^{-1}\cdot
L\|\bigr|\\
&=\sup_{L\in\mathcal L_n}
\bigl|\log_2\|L\|-\log_2\|\beta^{-1}\cdot L\|\bigr|\\
&=\widetilde c(\beta^{-1}).
\end{aligned}$$

For any $\beta_1,\beta_2\in B_n$ we have
$$\begin{aligned}
\widetilde c(\beta_1\beta_2)&=\sup_{L\in\mathcal L_n}
\bigl|\log_2\|\beta_1\beta_2\cdot L\|-\log_2\|L\|\bigr|\\
&=\sup_{L\in\mathcal L_n}
\bigl|\log_2\|\beta_1\cdot(\beta_2\cdot L)\|-\log_2\|\beta_2\cdot L\| \\
 & \hspace{1cm} +\log_2\|\beta_2\cdot L\|-\log_2\|L\|\bigr|\\
\end{aligned}
$$
$$
\begin{aligned}
&\leqslant\sup_{L\in\mathcal L_n}
\bigl|\log_2\|\beta_1\cdot L\|-\log_2\|L\|\bigr| \\
 & \hspace{1cm} +\sup_{L\in\mathcal L_n}\bigl|
\log_2\|\beta_2\cdot L\|-\log_2\|L\|\bigr|\\
&=\widetilde c(\beta_1)+\widetilde c(\beta_2).
\end{aligned}$$

Thus the formula
$$\dl(\beta_1,\beta_2)=\widetilde c(\beta_1^{-1}\beta_2)=
\sup_{L\in\mathcal L_n}\bigl|\log_2\|\beta_1^{-1}\cdot L\|-
\log_2\|\beta_2^{-1}\cdot L\|\bigr|$$
defines another metric $\dl$ on $B_n$.
The proof of the following claim is easy.

\begin{prop}
For either of the metrics $\dl$ and $\dw$, the standard injection
$B_n\rightarrow B_{n+1}$ is distance-preserving.
\end{prop}

The main result of this section is

\begin{theorem} \label{T:metricsqi}
The metrics $\dl$ and $\dw$ on $B_n$ are quasi-isometric,
namely, for any $\beta_1,\beta_2\in B_n$, $\beta_1\ne\beta_2$,
the following holds:
$$\log_32\leqslant\frac{\dw(\beta_1,\beta_2)}{\dl(\beta_1,\beta_2)}
\leqslant9n.$$
\end{theorem}

\begin{proof}
The first inequality is deduced from Lemma~\ref{L:one-Delta}
by analogy with the proof of the easy part of
Theorem~\ref{T:main1}.
The second inequality follows from Theorem~\ref{T:main2}.
\end{proof}

One can immediately deduce an analogue result for the group 
$B_n/\langle \Delta^2\rangle$, i.e., the quotient of the $n$ string
braid group by its centre. This is the finite index subgroup of the 
mapping class group of the $n+1$ times punctured sphere $S^2_{n+1}$
consisting of those elements which fix the $n+1$st puncture.

The metric $\dw$ on $B_n$ induces a metric on 
$B_n/\langle \Delta^2 \rangle$, which we shall still denote $\dw$.
So by definition the $\dw$-distance of an element 
$\beta$ of $B_n/\langle \Delta^2\rangle$ from the trivial element is
$\min_{k\in\Z} \dw(\beta \Delta^{2k},1_{B_n})$.
Similarly, we can define an analogue of the metric $\dl$ for 
$B_n/\langle \Delta^2\rangle$. Geometrically, this means that  two
laminations on $D_n$ are regarded as equivalent if one can be deformed into 
the other, where the  deformation must preserve $\partial D_n$ setwise, 
but not necessarily pointwise: all the $n-1$ endpoints of arcs can be slid
simultaneously through an angle of $2\pi z$ ($z\in\Z$) along the boundary. 
This modified notion of equivalence yields
a modified notion of complexity of a lamination, and thus an analogue
of the metric $\dl$ on $B_n/\langle \Delta^2 \rangle$. Again, this
metric shall still be denoted $\dl$. As an immediate consequence of
theorem \ref{T:metricsqi} we have

\begin{cor}\label{C:twometrics'}
The metrics $\dl$ and $\dw$ on $B_n/\langle \Delta^2\rangle$ are 
quasi-isometric.
\end{cor}


\subsection{\Tei spaces}\label{S:Teich}

Among the handful of metrics which are habitually imposed upon the
\Tei space $\calT(S)$ of a surface $S$, there are notably the 
\emph{\Tei metric} (which can be interpreted in terms of stretch factors 
of extremal lengths of curves on the surface, see \cite{Kerckhoff}), and 
W.~Thurston's \emph{Lipschitz metric} \cite{AmericanDiet}. The latter
metric can be interpreted in terms of stretch factors of 
hyperbolic lengths of curves on the surface. This interpretation, which 
is due to Thurson, will be recalled below. A theorem 
of Choi and Rafi \cite{ChoiRafi} states that the distance $d(\sigma,\tau)$ 
between two points $\sigma, \tau$ in \Tei space according to the two metrics
are the same up to an additive constant, provided the two points lie in
the \emph{thick part} of \Tei space, meaning that they possess no
hyperbolic geodesics shorter than the Margulis constant. In particular,
the thick parts of \Tei space, equipped with these two metrics, are
quasi-isometric.

It should be mentioned that the Lipschitz ``metric'' is not actually
symmetric, but could easily be turned into a metric by symmetrizing.
Moreover, the Lipschitz metric and its symmetrized version are
quasi-isometric on the thick part of \Tei space.

We recall that there is a natural action of the mapping class group
of $S^2_{n+1}$, and hence of $B_n/\langle \Delta^2\rangle$, on 
the \Tei space $\calT(S^2_{n+1})$, and this restricts to a cocompact 
action on the thick part $\calT_{\mathrm{thick}}(S^2_{n+1})$ of 
\Tei space. Thus for any fixed
point $\sigma_*$ in the thick part, the orbit of
$\sigma_*$ under the action provides an embedding
$\Phi\co B_n/\langle\Delta^2\rangle\to \calT_{\mathrm{thick}}(S^2_{n+1}), 
\beta\mapsto\beta\cdot\sigma_*$.
Let us now equip $\calT_{\mathrm{thick}}(S^2_{n+1})$ with the restriction
of the Teichm\"uller, or equivalently, the Lipschitz metric on the 
full \Tei space, see \cite{ChoiRafi}.
The aim of this section is to prove that the pullback metric on the
braid group is, up to quasi-isometry, either one of the metrics defined
in the previous section. Thus the metric spaces constructed in section
\ref{S:twometrics} turn out to be combinatorial models for the thick 
part of \Tei space.

\begin{prop}\label{P:Teich}
The embedding $\Phi \co (B_n/\langle\Delta^2\rangle,\dl) \to 
(\calT(S^2_{n+1}),d_{\mathrm{Lipschitz}})$
is quasi-isometric.
\end{prop}

It should be stressed that this result is quite easy to prove, and
similar results are already in the literature
(see e.g.~theorem 2.2 of \cite{ChoiRafi}). 
What is more surprising is that, using Corollary \ref{C:twometrics'}
and Choi and Rafi's comparison between \Tei and Lipschitz metric on 
\Tei space \cite{ChoiRafi,RafiCombModel}, we obtain

\begin{cor}
The following four spaces are mutually quasi-isometric:
$$ 
\begin{array}{ll}
(1) \ (B_n/\langle\Delta^2\rangle,\dw) & 
(2) \ (B_n/\langle\Delta^2\rangle,\dl) \\
(3) \ (\calT_{\mathrm{thick}}(S^2_{n+1}),d_{\mathrm{Lipschitz}}) & 
(4) \ (\calT_{\mathrm{thick}}(S^2_{n+1}),d_{\mathrm{Teichm.}})
\end{array}
$$
\end{cor}

The fact that repeated Dehn twists yield a logarithmically escaping path
in \Tei space can already be seen from \cite{Minsky}. For completeness,
we give a proof of proposition \ref{P:Teich}.

\begin{proof}[Proof of Proposition \ref{P:Teich}] 
We shall use the following notation. If $f,g\co X\to \R$
are two functions, where $X$ is any set, then we say $f$ and $g$
are \emph{comparable}, and write $f\comp g$, if there exist constants
$C\geq 1$ and $d\geq 0$ such that 
$\frac{1}{C}\cdot g(x)-d\leq f(x)\leq C\cdot g(x)+d$.

Now, for $\alpha$ an isotopy class of simple closed curves
in $D_n$, and $\sigma$ a hyperbolic structure on $D_n$ (i.e.,~a point
in $\calT$) we shall denote $l_\sigma(\alpha)$ the hyperbolic length
of the shortest representative of $\alpha$, measured in the metric $\sigma$.
According to Thurston \cite{AmericanDiet}, there are two equivalent
definitions of the Lipschitz metric, among them the following:
if $\sigma, \tau$ are two hyperbolic structures, then
$$d_{\rm Lipschitz}(\sigma,\tau)=
\sup_{\alpha} (\ \log(l_\sigma(\alpha))-\log(l_\tau(\alpha))\ )$$
where the supremum is taken over all simple closed curves in $D_n$.
In particular, if we are trying to measure the distance between
$\sigma$ and its translate under the action of a braid $\beta$, we obtain
$$
d_{\rm Lipschitz}(\sigma,\beta\cdot\sigma) = \sup_{\alpha}\, 
(\,\log(l_\sigma(\alpha))-\log(l_\sigma(\beta\cdot\alpha))\,)
$$
Now we recall the well-known fact that for any fixed point 
$\sigma_*\in\calT$, there exist constants $c,C>0$ such that 
$$
c\cdot l_{\sigma_*}(\alpha)\leq \|\alpha\| \leq C\cdot l_{\sigma_*}(\alpha)
$$
for all $\alpha$. That is, $l_{\sigma_*}(\alpha)$ and $\|\alpha\|$
are in bilipschitz correspondence, and in particular, we have 
$l_{\sigma_*}(\alpha)\comp \|\alpha\|$.
(The reason why this is true is that for \emph{simple} closed geodesic 
curves $\alpha$ in $S^2_{n+1}$, equipped
with the metric $\sigma_*$, there are lower and upper bounds for the
lengths of the components of intersection of $\alpha$ with the lower 
and upper half plane.)


Moreover, there are global lower bounds on $l_{\sigma_*}(\alpha)$
(namely the Margulis constant) and on $\|\alpha\|$ (namely $2$).
Thus we can deduce that the logarithms of these quantities are also
comparable:
$$
\log(l_{\sigma_*}(\alpha)) \comp \log(\|\alpha\|).
$$
Now let $\alpha_1,\ldots,\alpha_k$ denote any finite family of simple 
closed curves with the property that every simple closed curve in 
$S^2_{n-1}$, except those enclosing a single puncture, can be obtained 
from one of the $\alpha_i$s by the action of some braid. 
We calculate
\begin{eqnarray*}
d_{\rm Lipschitz}(\sigma_*,\beta\cdot\sigma_*) & \comp & \sup_{\alpha}\, 
\left(\,\log(\|\alpha\|)-\log(\|\beta\cdot\alpha\|)\,\right)\\
 & = & \sup_{\zeta\in B_n} \sup_{i=1,\ldots,k} 
\left(\ \log(\|\zeta\cdot\alpha_i\|)-
                           \log(\|\beta\zeta\cdot\alpha_i\|)\ \right)\\
 & \comp & \sup_{\zeta\in B_n} \left(\,\sum_{i=1}^{k} \left(
\log(\|\zeta\cdot\alpha_i\|)-
                     \log(\|\beta\zeta\cdot\alpha_i\|)\right)\,\right)\\
 & \comp & \sup_{\zeta\in B_n} \left( 
 \log\left(\sum_{i=1}^{k} \|\zeta\cdot\alpha_i\|\right)
 -\log\left(\sum_{i=1}^{k}\|\beta\zeta\cdot\alpha_i\|\right)\right).
\end{eqnarray*}

We shall fix one very particular choice for the
family $\alpha_1,\ldots,\alpha_k$ namely the one indicated in Figure 
\ref{F:alphai}---in particular, in our choice we have $k=n-1$.
\begin{figure}[htb]
\centerline{\input{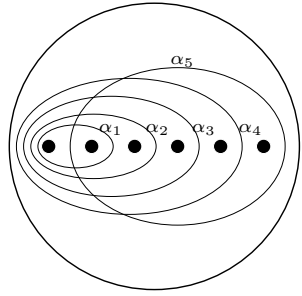}}
\caption{The curves $\alpha_1,\ldots,\alpha_{n-1}$}\label{F:alphai}
\end{figure}
With this particular choice for the family $\alpha_1,\ldots,\alpha_k$
we have the following comparison with the complexity of the curve 
diagram $\zeta\cdot E$:
$$
\|\zeta\cdot E\|   -n+1 
\leq \sum_{i=1}^{n-1} \|\zeta\cdot\alpha_i\| \leq \|\zeta\cdot E\|+n-1.
$$
In particular, we obtain
\begin{eqnarray*}
d_{\rm Lipschitz}(\sigma_*,\beta\cdot\sigma_*) & \comp & \sup_{\zeta\in B_n}
 \left( \log\left(\sum_{i=1}^{n-1} \|\zeta\cdot\alpha_i\|\right)
 -\log\left(\sum_{i=1}^{n-1}\|\beta\zeta\cdot\alpha_i\|\right)\right)\\
 & \comp & \sup_{\zeta\in B_n}
  \left( \log\left(\|\zeta\cdot E\|\right)
 -\log\left(\|\beta\zeta\cdot E\|\right)\right)\\
 & \comp & \dl(1_{B_n},\beta),
\end{eqnarray*}
which is what we wanted to prove.
\end{proof}


\subsection{Dehornoy braid ordering}\label{S:ordering}
In 1991 Patrick Dehornoy discovered that the braid group $B_n$
admits a left-invariant ordering \cite{Deh94}.
His purely algebraic construction was based on the following
notion of $\sigma$-positive braid word.

\begin{definition}
A braid word $w$ is said to be $\sigma_k$-positive (respectively,
negative), if it contains $\sigma_k$, but does not contain
$\sigma_k^{-1}$ and $\sigma_i^{\pm1}$ with $i<k$
(respectively, contains $\sigma_k^{-1}$, but
not $\sigma_k$ and $\sigma_i^{\pm1}$ with $i<k$).
If $w$ does not contain $\sigma_i^{\pm1}$ with $i\leqslant k$, it
is called $\sigma_k$-\emph{neutral}.

A braid word $w$ is
said to be $\sigma$-positive (respectively, $\sigma$-negative),
if it is $\sigma_k$-positive (respectively, $\sigma_k$-negative)
for some $k\leqslant n-1$. A braid word $w$ is said to be
$\sigma$-consistent if it is either trivial or $\sigma$-positive,
or $\sigma$-negative.
\end{definition}

\def\theenumi{{\rm\arabic{enumi})}}
\begin{theorem}[Dehornoy \cite{Deh94}]
For any braid $\beta\in B_n$, exactly one of the following is true:
\begin{enumerate}
\item
$\beta$ is trivial;
\item
$\beta$ can be presented by $\sigma_k$-positive braid word for some $k$;
\item
$\beta$ can be presented by $\sigma_k$-negative braid word for some $k$.
\end{enumerate}
In the latter two cases $k$ is unique.
\end{theorem}

Thus, it makes sense to speak about $\sigma$-positive
and $\sigma_k$-positive (or $\sigma$-, $\sigma_k$-negative) braids.
It is then an immediate consequence that the relation $<$ on $B_n$
defined by the rule: $\beta_1<\beta_2$ if and only if
$\beta_1^{-1}\beta_2$ is $\sigma$-positive, is a left-invariant
linear ordering on $B_n$, see~\cite{Deh94}.

It was noticed in~\cite{FGRRW} that this notion of $\sigma$-positivity
has a nice geometric interpretation in terms of curve diagrams.
We refer the reader to the monograph~\cite{DDRW} for a survey
of this and many other explanations of the phenomenon and
different proofs of Dehornoy's theorem.

Although many approaches to $\sigma$-ordering have been developed since
Dehornoy's discovery, the following question is currently unsettled:
is there a polynomial in $\ell$ which is an upper bound on the length
of the shortest $\sigma$-consistent braid word representing
a braid of length $\ell$?
Dehornoy's original algorithms (in \cite{Deh94},
and handle reduction \cite{Deh97}) and the algorithm from~\cite{FGRRW}
give only an exponential bound on the length of
the shortest $\sigma$-consistent representative.

At the end of the paper we shall present some further reasons for
believing that a linear bound exists.
The aim of the current section is to solve a closely related problem,  
namely, to give a positive answer to the question above with the 
ordinary braid length replaced by the $\Delta$-length.
This assumes the following extension of the notion
of $\sigma$-positive braid word: a word in the alphabet
$\{\Delta_{ij}\}_{0<i<j<n}$ is said to be $\sigma$-positive
if, for some $k<l$, it contains $\Delta_{kl}$, and
contains neither $\Delta_{kj}^{-1}$ nor $\Delta_{ij}^{\pm 1}$
with $i<k$ and any $j$. In other words, a word $w$ in letters
$\Delta_{ij}$ is $\sigma$-positive (negative, neutral) if the word
in standard generators $\sigma_i$ obtained from $w$ by
expansion~(\ref{E:delta}) is.

\begin{theorem}\label{T:sigma-consistency}
Any braid $\beta\in B_n$ can be presented by a $\sigma$-consistent
word $w$ in the alphabet $\{\Delta_{ij}\}$ such that
\begin{equation}
\ll(w)\leqslant30n\cdot\ll(\beta).
\end{equation}
\end{theorem}

The following lemma plays a key r\^ole in the proof.
Denote by $E_2$ the lamination in $D_n$ whose closure
is a circle in the complex plane
surrounding the leftmost puncture and
the leftmost point $*$ of $D_n$. (The notation
is motivated by the fact that this circle
coincides with the trivial curve diagram $E$ in
the case $n=2$.)

\begin{lemma}\label{L:sigma-consistency}
Let $\beta\in B_n$ be a $\sigma_1$-positive braid,
and let $u$ be the braid word spelt out by the untangling
procedure of Section~\ref{S:optimize} applied
to the lamination $\beta\cdot E_2$.
Then the braid word $u$ is $\sigma_1$-negative,
the braid $u\cdot\beta$ is $\sigma_1$-neutral,
and we have
\begin{equation}\label{E:u}
\ll(u)\leqslant\log_23\cdot9n\cdot\ll(\beta).
\end{equation}
\end{lemma}

\begin{proof}
First of all, we remark that the algorithm of Section~\ref{S:optimize}
constructs a sequence of puncture slidings rather than a concrete
braid word. In some cases such a sliding can be written in two different
ways, an example is shown in Fig.~\ref{F:sliding}. More precisely then,
what we are going to prove is that, \emph{under an appropriate
choice} of the decomposition of the slidings into a product of $\Delta$s
at each step of the untangling procedure, we get a $\sigma_1$-negative
word $u$.

A sliding of a group of punctures, like the one shown in
Fig.~\ref{F:sliding},
is not $\sigma_1$-neutral if and only if one of the following occurs:
the leftmost puncture is slid; or punctures are slid over or under
the leftmost (true) one. If none of these takes place, it does not matter
which decomposition, $\Delta_{ij}^\epsilon\Delta_{ij'}^{-\epsilon}$
or $\Delta_{ij}^\epsilon\Delta_{i'j}^{-\epsilon}$,
to choose---both are $\sigma_1$-neutral.

If a sliding of a group of punctures is not $\sigma_1$-neutral, then
its decomposition of the form $\Delta_{ij}^\epsilon\Delta_{i'j}^{-\epsilon}$
is $\sigma_1$-definite, and it is $\sigma_1$-positive (respectively,
negative) if and only if the punctures are slid
clockwise (respectively, counterclockwise).

Thus, in order to prove that $u$ is $\sigma_1$-negative it suffices
to show that
all the clockwise slidings that occur during untangling the lamination
$\beta\cdot E_2$ are $\sigma_1$-neutral,
but the whole word $u$ is not $\sigma_1$-neutral.

By construction, the braid $u\cdot\beta$ preserves $E_2$, which implies
that this braid is $\sigma_1$-neutral. Since $\beta$ is not
$\sigma_1$-neutral, the braid represented by $u$ cannot be
$\sigma_1$-neutral, either. Let us check that $\sigma_1$-positive
slidings do not appear in $u$.

To this end, we must examine all the clockwise slidings and make
sure that they do not involve the leftmost puncture, i.e.\
the arcs of all clockwise sliding are always on the right of the leftmost
puncture.

Similarly to~\cite{FGRRW} one can show that $\sigma_1$-positivity
of the braid $\beta$ is equivalent to the following property of the
lamination $L=\beta\cdot E_2$, which is assumed to be tight with
respect to $\R$:
$$
\parbox{0.8\textwidth}{there is an arc $\gamma\subset\overline L$ lying in
the lower half-plane with endpoints at $\R$ such that
the segment of $\R$ bounded by the endpoints of $\gamma$ contains $*$
and does not contain a true puncture.}
\eqno(\star)$$
See also~\cite{DDRW} for an explanation of the relation between
$\sigma$-positivity and a ($\star$)-like property of curve diagrams.

The lamination $\overline L$
intersects the ray $(-\infty,*)\subset\R$ exactly
once, at the left endpoint of $\gamma$. (This means, in particular,
that there is only one arc $\gamma$ satisfying~($\star$).)
This implies that no puncture is slid out of the disk $D_n$
during untangling $L$. Indeed, in the general case,
punctures can be slid out of the disk at the relaxation
step that follows an AA- or AB-transmission. One can see that,
in both situations, any arc in the upper half-plane along which
a puncture is slid counterclockwise must lie in between
two parallel arcs of the lamination.
Thus, sliding a puncture out of the disk may be forced only if
the closure of the lamination $L$ intersects the ray $(-\infty,*)$
at least twice, which is not the case for $L=\beta\cdot E_2$.

Thus, the braid $\beta_{p'+1}'$ that is composed of clockwise
slidings and is a potential source of a $\sigma_1$-positive contribution
to $u$ is trivial.

Further, we claim that property ($\star$) holds during the whole
untangling process. Indeed, the arc $\gamma$ cannot be
essential, so it always remains untouched. Thus property ($\star$)
is violated only when a true puncture is slid to a point in between the
endpoints of $\gamma$, which is the moment when the lamination
gets untangled completely.

Now we can show that property ($\star$) guarantees that
no clockwise sliding involves the leftmost puncture.
That is, we have to exhibit a true puncture to the left of all
arcs along which clockwise slidings occur.
This is done by revising, case by case, the relaxations following
different types of transmissions. The AB spiralling and non-spiralling
cases are trivial because the corresponding slidings are
counterclockwise. In all the other cases, just before the transmission,
there must be an AA-strip whose left base is further to the left
than the right base of the numerator.
Let $R$ be the innermost such strip
($R$ may be the numerator in the AA case).

According to our transmission-relaxation rules, there
must be a puncture (true or false) $P$, say, between the bases of $R$.
In each individual case it is easy to show that $P$
is not obstructing and that it
lies to the left of the arcs along which clockwise slidings occur.
Thus if $P$ is a true puncture, then the clockwise slidings
are $\sigma_1$-neutral. We conclude by noticing
that $P$ cannot be the false puncture because otherwise
property ($\star$) would imply the presence of a circle
surrounding just $*$, which is absurd.

It remains to prove estimation (\ref{E:u}). This is done
by applying the first inequality in Theorem~\ref{T:main1}
and Theorem~\ref{T:main2},
and using the fact that $u$ is precisely the braid word
for which the estimation from Theorem~\ref{T:main2} has been obtained:
$$\ll(u)\leqslant9n\cdot\log_2\|L\|
\leqslant9n\cdot(\log_23\cdot\ll(\beta)+\log_2\|E_2\|)=
9n\cdot\log_23\cdot\ll(\beta),$$
since $\|E_2\|=1$.
\end{proof}

\begin{lemma}\label{L:sigma-free}
If a braid $\beta$ is $\sigma_1$-neutral then the shortest
braid word representing $\beta$ (where either the ordinary
braid length or the $\Delta$-length is used) is $\sigma_1$-neutral.
\end{lemma}

\begin{proof} If $\beta$ is $\sigma_1$-neutral, then the string that
starts in the leftmost position also ends in the leftmost position.
Now if $w$ is any braid word representing $\beta$, then we can obtain
another braid word $\widetilde{w}$ representing $\beta$ which is
$\sigma_1$-neutral and at most as long as $w$ in the following way:
we delete the string starting and ending in the leftmost position,
and replace it by a string that lies entirely to the left of all the
other (true) strings. The word $\widetilde w$ is then obtained
from $w$ by removing some letters $\sigma_i^{\pm1}$,
shifting indices of others:
$\sigma_i^{\pm1}\mapsto\sigma_{i-1}^{\pm1}$,
$\Delta_{ij}^p\mapsto\Delta_{i-1,j-1}^p$,
and replacing some $\Delta_{ij}^p$ with $\Delta_{i,j-1}^p$.
If $w$ was not $\sigma_1$-neutral, then $\widetilde w$
is strictly shorter than $w$.
\end{proof}

\begin{proof}[Proof of Theorem~\ref{T:sigma-consistency}]
Without loss of generality we may assume that
$\beta$ is $\sigma_k$-positive with some $k<n$.
The proof is by induction on $k$.
We skip the details of the induction step, which is an
easy consequence of Lemma~\ref{L:sigma-free},
and show that the assertion of the theorem
holds for a $\sigma_1$-positive braid.

Indeed, take $u$ from Lemma~\ref{L:sigma-consistency}.
Since the braid $\beta'=u\cdot\beta$ is $\sigma_1$-neutral,
Lemma~\ref{L:sigma-free} imlpies that $\beta'$
can be presented by a $\sigma_1$-neutral braid word
$v$ of $\Delta$-length
$$\ll(v)\leqslant\ll(u\cdot\beta)\leqslant\ll(u)+\ll(\beta)\leqslant
(\log_23\cdot9n+1)\cdot\ll(\beta).$$
The braid word $w=u^{-1}v$ is then $\sigma_1$-positive,
it represents $\beta$ and has $\Delta$-length
$$\ll(w)\leqslant\ll(u)+\ll(v)\leqslant(2\log_23\cdot9n+1)\cdot
\ll(\beta)\leqslant30n\cdot\ll(\beta).$$
\end{proof}


\subsection{Algorithmical issues}\label{S:algorithm}

The proof of Theorem~\ref{T:main2} presented above consists
in an explicit description of an algorithm that,
given an integral lamination $L$, finds a word-representative
of a braid that untangles $L$. In this section we
discuss the efficiency of the algorithm and, more generally,
of the algorithmical treatment of the braid groups based on it.

First of all, we remark that our estimations for the running time
of algorithms will be made for the computational model called
\emph{Random Access Memory Machine}. Roughly speaking, this
means that we assume the input to be in a reasonable range,
and estimate the number of elementary operations of a realistic
computer needed to implement the algorithm.
More precesily, we assume that the number of strands
$n$ is ``small enough'', so that its record fits one standard unit
of memory, and that any arithmetic operation
on integers between $0$ and $n$ takes constant time.
This is a reasonable assumption because actually it allows $n$
to be very large: if, say, four bytes are used to store an integer
(which is quite usual), then $n$ can be as large as $2^{31}$.
For other integers (which are not indices) we will assume
that their logarithm is ``small'' (i.e.\ smaller than $2^{31}$).

Since even for reasonably long braid words the implementation
of our algorithms may need to operate with ``large'' integers,
we will pay attention to
the number of elementary machine operations that are needed to
perform an arithmetic operation on those integers.
The most frequently used operations will be addition, subtraction,
and comparison. They consume logarithmic time
in the value of the larger operand, and we call them \emph{simple}.
Sometimes we will also need to perform
divisions and multiplications. We use the fact that dividing $k$ by $l$
takes $O(\log l\cdot\log(k/l))$
elementary operations, whereas multiplying $k$ by $l$ consumes
$O(\log k\cdot\log l)$ operations.

The next important question is how the input/output data and the objects
used in the algorithm are presented. Our main objects are: braid
words, laminations, and strip systems, so we briefly discuss
their numerical presentations.

We will assume that words of the form~(\ref{Delta-form})
are presented by the corresponding sequences of integers:
$$s;i_1,j_1,k_1,i_2,j_2,k_2,\dots,i_s,j_s,k_s.$$
Here $s$ and $i_t,j_t$ with $t=1,\dots, s$ are ``small'' integers,
whereas $k_t$, $t=1,\dots,s$, can be ``large''.
This implies that the length of such presentation of $w$
is of order $O(\ll(w))$.

The set $\mathcal L_n$ of laminations in $D_n$ can be naturally
identified with $\Z^{2n-2}\setminus\{0\}$ so that
the norm $\|L\|$ will become a norm in $\Z^{2n-2}$, and the action of each
generator $\sigma_i^{\pm1}\in B_n$ will be given by a finite
number of simple arithmetic operations on the coordinates of
the lamination. See~\cite[Chapter 8]{DDRW} and \cite{ybm} for details,
where a slightly different definition of laminations is used,
which results in two additional coordinates appearing in the ``code''
of a lamination. (In order for the formulas in~\cite{ybm,DDRW}
to work in our current settings, one should set the
two additional coordinates to $a_n=0$, $b_n=+\infty$.)
This implies the following

\begin{prop}
There exists an algorithm $A_1$ that, given a word $w$
in the generators $\Delta_{ij}$ representing a braid $\beta$,
computes the curve diagram
$\beta\cdot E\in\mathcal L_n=\Z^{2n-2}\setminus\{0\}$
in time $O(\ell(w)\cdot\ll(w)+n)$.
\end{prop}

The algorithm $A_1$ expands the given word $w$ by using~(\ref{E:delta}),
thus obtaining a word $w'$ in $\sigma_i^{\pm1}$-generators
recorded in the usual way. Then it generates the initial
lamination $E$ and applies, one after another, the letters of
$w'$ (from right to left) to the lamination.

This may be very inefficient if $w$ contains a subword $\Delta_{ij}^N$
with a very large $N$.
However, the action of the braid $\Delta_{ij}^N$ on a lamination $L$
can be computed without expanding the braid into
a product of $\sigma_i$s.
\begin{lemma}\label{L:DeltaN}
The action of $\Delta_{ij}^N$ on $L$ can be computed
in $O((n+\log N)\cdot\log\|L\|)$ operations.
\end{lemma}

\begin{proof}
Let us look at
the sequence of laminations $\Delta_{ij}^k\cdot L$, where
$k=\dots,-2,-1,0$, $1,2,\dots$. For large $|k|$ the laminations
$\Delta_{ij}^k\cdot L$ have a big ``spiral'' surrounding punctures
$i$ through $j$, and the ``thickness'' of this spiral grows linearly
with $k$. More formally, this means the following.

Let $k_0$ be an integer for which $\|\Delta_{ij}^{k_0}\cdot L\|$
is as small as possible. Let us cut $L$ along the real axis and
count the number of the obtained arcs having one endpoint between
the $i$th and $j$th punctures, and the other endpoint outside
this segment. In a sense, this is twice the number of strings
involved in the spirals of $\Delta_{ij}^N\cdot L$.
Let this number be $m$ and let $L_0$
be the lamination consisting of $m$ circles surrounding punctures
$i$ through $j$. Let us think of laminations as points in
$\Z^{2n-2}$. Then for any $p\geqslant1$ the following holds:
$$\Delta_{ij}^{k_0+2p}\cdot L=\Delta_{ij}^{k_0+2}\cdot L+(p-1)\cdot L_0,
\qquad
\Delta_{ij}^{k_0-2p}\cdot L=\Delta_{ij}^{k_0-2}\cdot L-(p-1)\cdot L_0.$$

So, we start by establishing the structure of a spiral in $L$
surrounding punctures $i$ through $j$, if there is one. Even if there
is no spiral, we compute $m$, i.e.\ the lamination $L_0$.
If there is a spiral, we also need to find its ``thickness'' $\theta$
and its direction (clockwise or counterclockwise).

By using flips of triangulations
in a similar way as described in~\cite{DDRW}, we can do all this
job in $O(j-i)\leqslant O(n)$ simple operations on integers
of order $\|L\|$. So, the structure of the spiral can be discovered
for $O(n\cdot\log\|L\|)$ elementary operations.

Depending on the direction of the spiral and the sign of $N$
different cases are possible.
It may happen that $\Delta_{ij}^N$ twists the spiral further,
in which case we are lucky, because we have $\Delta_{ij}^N\cdot L=
L+N/2\cdot L_0$, provided that $N$ is even. If $N$ is odd,
we shall also need to apply one $\Delta_{ij}$ ``explicitly'',
which, by the same ``flip argument'', takes $O(n)$ simple operations
on integers of order $\|L\|$. So, the total work in this case
is $O(n\cdot\log\|L\|)+O(\log N\cdot\log\|L\|)$, where
the second summand appears because we need to multiply $L_0$ by $N/2$.
(Note that $L_0$ has only two non-zero coordinates.)

The same estimation works if $\Delta_{ij}^N$ untwists the spiral
partially. The most involved case is when $\Delta_{ij}^N$ untwists
the spiral completely and then twists in the opposite direction.
In addition to the previous cases, we shall need to apply a few
more $\Delta_{ij}$s explicitly, and compute the number of
twists in the original spiral. The latter is done by computing
$[\theta/m]$, which consumes $O(\log N\cdot\log\|L\|)$ elementary
operations (because $\theta/m<N$).\end{proof}

By using induction we deduce the following from Lemma~\ref{L:DeltaN}.

\begin{prop}
There exists an algorithm $A_2$ that computes the curve diagram
$\beta\cdot E$ of a braid $\beta$ given by a braid word $w$ in time
$O(n\cdot\ll(w)^2)$.
\end{prop}

In order to implement the algorithm of section~\ref{S:optimize}
one needs to choose a presentation method for strip systems.
The most straightforward way to present a strip system $(L,S)$
is to provide coordinates of $L$, list all interval identifications
of $S$, and specify the positions of the punctures.
However, in order to make the algorithm more efficient it
is useful to include even more information in the object.
For example, one may keep a bi-directed
list of ``significant'' points of the axis, which are positions of
punctures and the endpoints of the bases of strips, and a collection
of cross-references between those points and the related objects
(punctures, bases of strips), so as to be able, say, for any base
of a strip to find the ``next'' one in a bounded number
of simple operations. We skip the boring details.

One can show that for an appropriate encoding of strip systems,
each non-spiralling transmission and the subsequent relaxation
described in Sections \ref{S:AHT}--\ref{S:optimize}
can be performed in $O(n)$ simple arithmetic operations whose
operands are of order $O(\|L\|)$. For performing a $d$-times
spiralling transmision on a strip system $(L,S)$ we additionally
need to implement one division (the width of the numerator is devided
by the sum of widths of the bases of denominators participating in
the transmission), which consumes $O(\log d\cdot\log\|L\|)$ operations.
Together with Theorem~\ref{T:main2} this implies the following.

\begin{prop}
There exist algorithms $A_3,A_4$ such that
\begin{enumerate}
\item
given a lamination $L$, $A_3$
computes a braid untangling $L$ in time $O(n^2\cdot(\log\|L\|)^2)$,
thus detecting whether $L$ is the curve diagram of some braid;
\item
given the curve diagrams $L_1$, $L_2$ of braids $\beta_1$, $\beta_2$,
$A_4$ computes the curve diagram of $\beta_1\circ\beta_2$
in time $O(n^2\cdot(\ll(\beta_1)+\ll(\beta_2))\cdot\ll(\beta_1))$.
\end{enumerate}
\end{prop}


We shall finish this paper with some remarks concerning open problems 
and possible further developments of our results.

Firstly, all our results concern punctured disks and spheres and
their mapping class groups. It would be useful to find generalizations 
applying to mapping class groups of more general surfaces.

Secondly, we conjecture that the untangling procedures defined in
sections \ref{S:relax} and \ref{S:ordering} describe paths in the
Cayley graph of $B_n$ which are uniform quasigeodesics with respect
to the standard metric (not our $\Delta$-metric) on $B_n$. Indeed, 
these paths look very much like train track splitting sequences, 
which are known to be quasigeodesics by a theorem of Hamenst\"adt 
(\cite{Ham} Proposition 3.1). However, the exact technical conditions
of Hamenst\"adt's theorem, and in particular the genericity condition, 
are not easy to satisfy. Our conjecture would in particular imply that 
every braid has a $\sigma_1$-consistent representative whose length is
bounded linearly by the length of the braid---the existence of such
a representative is still an open problem \cite{Deh97,DDRW}. If the
conjecture were true, then our untangling paths would have the interesting
property that they are short with respect to both the usual, and the
$\Delta$-metric on $B_n$.


Thirdly, it might be useful to give substance to the intuition that
every ``spiral'' that appears during our untangling algorithm is 
somehow ``visible'' in every reasonably short representative of the
braid, and in particular in the Garside normal form. The idea here
is that spirals correspond to passages of the \Tei geodesic through
the thin part of Teichm\"uller space.

Finally, there might be applications of our results to the study of sets of 
``short'' elements in the conjugacy class of a braid---for instance, 
the super summit set of a braid. The reason for this hope
is that conjugacy classes of braids correspond to free homotopy
classes of closed curves in moduli space.


\subsection*{Acknowledgments.} We thank Ian Agol, who, in a conversation with
I.\ Dynnikov, first suggested applying the techniques from \cite{AHT}
to our problem. Jason Behrstock suggested the connection with \Tei
spaces, and made some very helpful remarks on that subject.
We also thank Lee Mosher for very helpful discussions.
The work of I.~Dynnikov was supported in part by
Russian Foundation for Basic Research (grant no.~02-01-00659).
Finally, we thank the CNRS and the Russian Academy of Sciences
for their financial support: their joint Franco-Russian exchange program
paid for a two-week stay of B.~Wiest at the Steklov Institute, Moscow, and
for a one-week stay of I.~Dynnikov at Rennes University.

\end{document}